\documentclass{article}
\usepackage{amsmath,amsthm}
\usepackage{amssymb}
\usepackage{eepic}
\usepackage{epsfig}
\usepackage{array}
\usepackage{picinpar}
\usepackage{graphics}
\usepackage[mathscr]{eucal}
\usepackage{multicol,multirow}

\allinethickness{0.4mm}

\theoremstyle{plain}
\newtheorem{theorem}{Theorem}
\newtheorem{lemma}{Lemma}
\newtheorem{prop}{Proposition}
\newtheorem{corollary}{Corollary}

\newcommand{\Z}{\mathbb Z}
\newcommand{\Stg}{\mathop{\mathsf{St}}\nolimits_G}
\newcommand\s{{\sigma}}
\newcommand\St{{\mathsf{St}}}
\newcommand{\Aut}{\mathop{\rm Aut}\nolimits}
\newcommand{\Isom}{\mathop{\rm Isom}\nolimits}

\newcommand\Sym{{\mathsf{Sym}}}
\newcommand{\lims}[1][G]{\mathscr{J}_{#1}}
\renewcommand\H{{\mathcal H}}

\title{Groups generated by $3$-state automata
over a $2$-letter alphabet, II}
\author{Ievgen Bondarenko, Rostislav Grigorchuk, Rostyslav Kravchenko,\\ Yevgen Muntyan, Volodymyr
Nekrashevych,\\ Dmytro Savchuk and Zoran \v{S}uni\'{c}\thanks{All
authors were partially supported by at least one of the NSF grants
DMS-0456185, DMS-0600975, and DMS-0605019} }
\date{\today}

\begin{document}

\maketitle

\begin{abstract}
Classification of groups generated by 3-state automata over a
2-letter alphabet started in \cite{bondarenko_gkmnss:clas32} is
continued.
\end{abstract}

\section{Introduction}

A comprehensive effort towards classification of all groups
generated by invertible automata on three states acting on 2 letters
($(3,2)$-automata) was started by our research group at Texas A\&M
University in~\cite{bondarenko_gkmnss:clas32} and the current text
presents further results in this direction.

The study of groups generated by finite automata is not new (it
started in the beginning of 1960's). At the moment, the subject is
going through a rather productive and mature phase in the sense that
it is both contributing solutions to outstanding problems in various
branches of mathematics (group theory, amenability, complex
dynamics, etc.) and is developing and discovering its own
fundamental concepts and theory leading to its own deep problems. A
short account of the history of ideas related to the subject, as
well as the most significant and interesting achievements of the
theory can be found in~\cite{bondarenko_gkmnss:clas32}. We also
recommend the article~\cite{gns00:automata} and the
book~\cite{nekrash:self-similar} to readers interested in becoming
more familiar with automaton groups.

Groups defined by automata on 2 states acting on two letters are
already classified.

\begin{theorem}[\cite{gns00:automata}]\label{thm:class22}
There are, up to isomorphism, $6$ different groups generated by
$2$-state automata over a $2$-letter alphabet. Namely, the trivial
group, the cyclic group $C_2$ of order 2, Klein's Viergruppe
$C_2\times C_2$, the infinite cyclic group $\Z$, the infinite
dihedral group $D_\infty$ and the lamplighter group $\Z \wr C_2$.
\end{theorem}

We note that there are $2^3 \times 3^6 = 5832$ different (labeled)
automata on three states acting on two letters. Obvious symmetries
such as permutation of states, permutation of letters, and inversion
of states, together with minimization, reduces the number of
automata that needs to be checked to 194
(see~\cite{bondarenko_gkmnss:clas32} for more details). Through
additional work we were able to reduce the number of non-isomorphic
groups defined by these automata.

\begin{theorem}
There are no more than $124$ pairwise non-isomorphic groups defined
by automata on 3 states over a 2-letter alphabet.
\end{theorem}

For this text we chose 15 automata and we provide the basic
information on each group defined by these automata.
Section~\ref{groups} contain a list of these 15 automata along with
the basic information in a tabular form. Section~\ref{proofs}
provides proofs of the claims in the tables as well as other
information.

There are very few results that hold for the whole class of groups
defined by automata on 3 states over a 2-letter alphabet. The reason
is that this class contains many groups of rather distinct nature.
The complete list contains finite groups, virtually free abelian
groups, solvable groups (such as the lamplighter group $\Z \wr C_2$
and Baumslag-Solitar groups $BS(1,\pm3)$), the free group $F_3$ of
rank 3, the free product $C_2 \ast C_2 \ast C_2$, an amenable but
not sub-exponentially amenable group (the Basilica group), and some
examples about which we know very little.

\begin{theorem}
There are $6$ finite groups in the class: $\{1\}$ [1], $C_2$ [1090],
$C_2\times C_2$ [730], $D_4$ [847], $C_2\times C_2\times C_2$ [802]
and $D_4\times C_2$ [748].
\end{theorem}

\begin{theorem}
There are $6$ abelian groups in the class: $\{1\}$ [1], $C_2$
[1090], $C_2\times C_2$ [730], $C_2\times C_2\times C_2$ [802],
$\mathbb Z$ [731] and $\mathbb Z^2$ [771].
\end{theorem}

In addition, there are virtually free-abelian groups in the class
(having $\Z$, $\Z^2$ [2212], $\Z^3$ [752] or $\Z^5$ [968] as
subgroups of finite index).

\begin{theorem}
\label{thm:class_free} The only free non-abelian group in the class
is the free group of rank 3 generated by the Aleshin-Vorobets
automaton [2240]. Moreover, the isomorphism class of this automaton
group coincides with its equivalence class under symmetry.
\end{theorem}

For the precise definition of symmetric automata
see~\cite{bondarenko_gkmnss:clas32}.

\begin{theorem}
There are no infinite torsion groups in the class.
\end{theorem}


\section{Regular rooted tree automorphisms and self-similarity}

Fix a finite alphabet $X=\{0,\dots,d-1\}$ on $d$ letters, $d \geq
2$. The set of all words $X^*$ over $X$ has the structure of a
\emph{regular rooted $d$-ary tree}, which we also denote $X^*$. The
empty word is the \emph{root} and each vertex $w$ has $d$
\emph{children}, namely the words $wx$, for $x \in X$. \emph{Level}
$n$ of the tree $X^*$ is the set $X^n$ of all words of length $n$
over $X$.

The \emph{boundary} of the tree $X^*$, denoted $X^\omega$, is the
set of right infinite words over $X^*$. It corresponds to the set of
infinite geodesic rays in $X^*$ starting at the root.

Denote by $\Aut(X^*)$ the group of automorphisms of $X^*$. Note that
any tree automorphism fixes the root and preserves the levels of the
tree. Every automorphism $f$ of $X^*$ can be decomposed as
\begin{equation}\label{decomposition}
  f = \alpha_f(f_0,\dots,f_{d-1})
\end{equation}
where $f_x \in \Aut(X^*)$, for $x \in X$, and $\alpha_f$ is a
permutation of $X$. The action of $f$ on $X^*$ can be described in
terms of the decomposition~\eqref{decomposition} as follows. Each of
the $d$ automorphisms $f_x$, $x \in X$, in the $d$-tuple
$(f_0,\dots,f_{d-1})$ acts on the corresponding subtree $xX^*$
consisting of the words over $X$ starting in $x$, after which the
permutation $\alpha_f$ permutes these subtrees. In other words
\begin{equation} \label{fxw}
 f(xw) = \alpha_f(x) f_x(w),
\end{equation}
for a letter $x$ in $X$ and word $w$ over $X$. The automorphism
$f_x$, $x \in X$, are called the (first level) \emph{sections} of
$f$ (and are sometimes denoted by $f|_x$), while $\alpha_f$ is
called the \emph{root permutation} of $f$. The notion of a section
can be extended recursively as follows. The section of $f$ at the
root is $f$ itself. For a word $wx$, where $w$ is a word and $x$ a
letter over $X$, the section of $f$ at $wx$ is defined as the first
level section $(f_w)_x$ of $f_w$ at $x$.

The group $\Aut(X^*)$ decomposes algebraically as
\begin{equation}\label{eqn:iter_wreath}
 \Aut(X^*) = \Sym(X) \ltimes (\Aut(X^*) \times \dots \times \Aut(X^*)) = \Sym(X)\wr\Aut(X^*),
\end{equation}
where $\wr$ is the \emph{permutational wreath product} (the
coordinates of $\Aut(X^*) \times \dots \times \Aut(X^*) =
\Aut(X^*)^X$ are permuted by $\Sym(X)$). The product of two
automorphisms $f$ and $g$ of $X^*$ is given by
\[
 \alpha_f(f_0,\dots,f_{d-1}) \alpha_g(g_0,\dots,g_{d-1}) =
 \alpha_f\alpha_g (f_{g(0)}g_0,\dots,f_{g(d-1)}g_{d-1}).
\]

A group $G$ of tree automorphisms is \emph{self-similar} if, for
every $g$ in $G$ and a letter $x$ in $X$ there exists a letter $y$
in $X$ and an element $h$ in $G$ such that
\[ g(xw) = yh(w), \]
for all words $w$ over $X$. More simply put, a group $G$ of tree
automorphisms is self-similar if every section of every element in
$G$ is again an element of $G$.

The $n$-th \emph{level stabilizer} $\St_{G}(n)$ of a group $G \leq
\Aut(X^*)$ is the group of automorphisms in $G$ that fix (pointwise)
the vertices at level $n$ (which implies that the vertices at the
lower levels are also fixed).

A self-similar group $G$ is called \emph{self-replicating} if, for
every vertex $u$, the homomorphism $\varphi_u:\Stg(u) \to G$ from
the stabilizer $u$ in $G$ to $G$, given by $\varphi(g)=g_u$, is
surjective.

A group $G$ of tree automorphisms of $X^*$ is \emph{level
transitive} if it acts transitively on each level of the tree $X^*$.

A group $G$ of tree automorphisms is contracting if there exists a
finite set $\mathcal N\subset G$, such that for every $g\in G$,
there exists $N>0$, such that $g_v\in\mathcal N$ for all vertices
$v\in X^*$ of length at least $N$ (see~\cite{nekrash:self-similar}).
The minimal set $\mathcal N$ with this property is called the
\emph{nucleus} of $G$. For a finitely generated group $G=\langle S
\rangle$, where $S$ is finite, the contracting condition is
equivalent to the existence of constants $\kappa$, $C$, and $N$,
with $0 \leq \kappa < 1$, such that $\ell(g_v) \leq \kappa \ell(g) +
C$, for all vertices $v$ of length at least $N$ and $g \in G$ (the
length $\ell(g)$ refers to the word length of the element $g$ in $G$
with respect to $S$).

A group $G$ of tree automorphisms of $X^*$ is a \emph{regular weakly
branch group} over its subgroup $K$ if $G$ acts spherically
transitively on $X^*$, $K$ is normal subgroup of $G$, and $K \times
\dots \times K$ is \emph{geometrically contained} in $K$. The latter
means that the first level stabilizer $\St_K(1)$ contains a subgroup
of the form $K \times \dots \times K$, where each factor acts on the
corresponding subtree hanging below the first level of the tree. For
more on (weakly) branch groups
see~\cite{grigorch:branch,bar_gs:branch}.

\section{Limit spaces and Schreier graphs}

Let $G$ be a finitely generated self-similar contracting group of
tree automorphisms of $X^*$. Denote by $X^{-\omega}$ the space of
left infinite words over $X$. Two elements $\ldots x_3x_2x_1, \ldots
y_3y_2y_1\in X^{-\omega}$ are \emph{asymptotically equivalent} with
respect to the action of $G$, if there exist a finite set $F \subset
G$ and a sequence $\{f_k\}_{k=1}^\infty$ of elements in $F$ such
that
\[
f_k(x_kx_{k-1}\ldots x_2x_1)=y_ky_{k-1}\ldots y_2y_1
\]
for every $k\geq 1$.

The quotient space $\lims$ of the topological space $X^{-\omega}$ by
the asymptotic equivalence relation is called the \emph{limit space}
of the self-similar action of $G$.

The limit space $\lims$ is metrizable and finite-dimensional (its
dimension is bounded above by the size of the nucleus of $G$). If
$G$ acts spherically transitively, then the limit space $\lims$ is
connected.

We recall now the definition of Schreier graph. Let $G$ be a group
generated by a finite set $S$ and let $G$ acts on a set $Y$. The
\emph{Schreier graph} of the action $(G,Y)$ is the graph
$\Gamma(G,S,Y)$ whose set of vertices is $Y$ and set of edges is
$S\times Y$, where the initial vertex of the arrow $(s,y)$ is $y$
and its terminal vertex is $s(y)$. For $y\in Y$, the Schreier graph
$\Gamma(G,S,y)$ of the action of $G$ on the $G$-orbit of $y$ is
called the \emph{orbital Schreier graph} of $G$ at $y$.

Let $G=\langle S \rangle$, where $S$ is finite, be a finitely
generated subgroup of $\Aut(X^*)$ and consider the Schreier graphs
$\Gamma_n(G,S)=\Gamma(G,S,X^n)$. Let $\xi=x_1x_2x_3\ldots\in
X^{\omega}$. Then the pointed Schreier graphs
$(\Gamma_n(G,S),x_1x_2\ldots x_n)$ converge in the local topology to
the pointed orbital Schreier graph $(\Gamma(G,S,\omega),\xi)$
(see~\cite{grigorchuk-z:cortona} for more details on the topology of
the space of pointed graphs).

Schreier graphs are related to the computation of the spectrum of
the Markov operator $M$ on the group. Given a finitely generated
group $G=\langle S \rangle$, where $S$ is finite, acting on the tree
$X^*$, there is a unitary representation of $G$ in the space of
bounded linear operators $\H=B(L_2(X^\omega))$, given by
$\pi_g(f)(x)=f(g^{-1}x)$. The Markov operator $M$ associated to this
unitary representation is given by
\[
 M = \frac{1}{|S \cup S^{-1}|} \sum_{s \in S \cup S^{-1}} \pi_s.
\]
The spectrum of $M$ for a self-similar group $G$ can be approximated
by the spectra of the finite dimensional Markov operators $M_n$, $n
\geq 0$, related to the permutational representations of $G$
provided by the action on the levels of the tree. The union of the
spectra of $M_n$ approximates the spectrum of $M$ in the sense that
\[sp(M)=\overline{\bigcup_{n\geq0}sp(M_n)}.\]
For more on this see~\cite{bartholdi_g:spectrum}.

\section{Definition of automaton groups}\label{definition}

We now formally describe the way finite automata define finitely
generated self-similar groups of tree automorphisms.

A \emph{finite invertible automaton} $A$ is a quadruple
$A=(Q,X,\tau,\rho)$ where $Q$ is a finite set of \emph{states}, $X$
is a finite \emph{alphabet} of cardinality $d \geq 2$, $\tau:Q
\times X \to Q$ is a map, called \emph{transition map}, and $\rho:Q
\times X \to X$ is a map, called \emph{output map}, such that, for
each state $q$ in $Q$, the restriction $\rho_q: X \to X$ given by
$\rho_q(x)=\rho(q,x)$ is a permutation of $X$.

Each state $q$ of the automaton $A$ defines a tree automorphism of
$X^*$, also denoted by $q$, by declaring that the root permutation
of $q$ is $\rho_q$ and the section of $q$ at $x$ to be $\tau(q,x)$.
Therefore
\begin{equation}
\label{eqn:autom_action}
q(xw) = \rho_q(x) \tau(q,x)(w)
\end{equation}
for all states $q$ in $Q$, letters $x \in X$ and words $w$ over $X$.

The group of tree automorphisms generated by the states of an
invertible automaton $A=(Q,X,\tau,\rho)$ is called the
\emph{automaton group} defined by $A$ and denoted by $G(A)$.


The boundary $X^\omega$ of the tree $X^*$ is endowed with a natural
metric (infinite words are close if they have long common
beginnings). The group of isometries $\Isom(X^\omega)$ with respect
to this metric is canonically isomorphic to $\Aut(X^*)$. Therefore
the action of the automaton group $G(A)$ on $X^*$ can be extended to
an isometric action on $X^\omega$ using relation
\eqref{eqn:autom_action}, which holds for infinite words $w$ as
well.

An invertible automaton $A$ can be represented by a labeled directed
graph, called Moore diagram, in which the vertices are the states of
the automaton, each state $q$ is labeled by its own root permutation
$\rho_q$ and, for each pair $(q,x) \in Q \times X$, there is an edge
from $q$ to $q_x=\tau(q,x)$ labeled by $x$. Moore diagrams of the 15
automata on 3 states that are subject of this article are provided
in Section~\ref{groups}.

\section{Selected groups}\label{groups}

Information on selected groups generated by $(3,2)$-automata is
provided in this section. The list appearing here supplements the
list provided in~\cite{bondarenko_gkmnss:clas32} and we keep the
same notation.

\begin{itemize}
\item
Rels - a list of some relators in the group. In particular, lists of
independent relators of length up to $20$ are included. In many
cases, the given relations are not sufficient (some of the groups
are not finitely presented).
\item
SF - the size of the factors $G/\Stg(n)$, for $n\geq0$.
\item Gr - the values of the growth
function $\gamma_G(n)$, for $n\geq0$, with respect to the generating
system consisting of $a$, $b$, and $c$.
\end{itemize}

In addition, in each case, a histogram for the spectral density of
the operator $M_9$ corresponding to the action on level $9$ of the
tree is provided for every automaton. Approximation of the limit
space is provided for the automaton 775.

The relations (and some other data) were obtained using the computer
algebra system GAP~\cite{GAP4} and software developed by Y.~Muntyan
and D.~Savchuk. \vspace{0.5cm}

\noindent\begin{tabular}{@{}p{172pt}p{174pt}@{}}
\vspace{0.1cm}

\textbf{Automaton number $741$} \vspace{.2cm}

\begin{tabular}{@{}p{48pt}p{200pt}}

$a=\sigma(c,a)$

$b=(b,a)$

$c=(a,a)$& Group:

Contracting: \textit{no}

Self-replicating: \textit{yes}\\
\end{tabular}

\begin{minipage}{230pt}
\setlength{\fboxrule}{0cm}
\begin{window}[5,r,{\fbox{\shortstack{\hspace{1.6cm}~\phantom{a}\\ \vspace{3.8cm}\\ \phantom{a}}}},{}]
Rels: $ca^{2}$, $b^{-1}a^{-1}cb^{-1}ababa$,\\
$bab^{-1}ca^{-1}b^{-1}aba$, $a^{-1}b^{-1}a^{-1}b^{-1}cabc^{-1}ab$,\\
$a^{-1}b^{-1}a^{-1}b^{-1}acbc^{-1}ab$, $b^{-1}c^{-3}b^{-1}cbcbc$,\\
$a^{-1}b^{-1}abc^{-1}a^{-1}bab^{-1}c$,
$a^{-1}b^{-1}aba^{-1}c^{-1}bab^{-1}c$.\\
\\
SF: $2^0$,$2^{1}$,$2^{3}$,$2^{6}$,$2^{12}$,$2^{23}$,$2^{45}$,$2^{88}$,$2^{174}$\\
Gr: 1,7,29,115,441,1643\\
\end{window}
\end{minipage}
& \hfill~

\hfill
\begin{picture}(1450,1090)(0,130)
\put(200,200){\circle{200}} 
\put(1200,200){\circle{200}}
\put(700,1070){\circle{200}}
\allinethickness{0.2mm} \put(45,280){$a$} \put(1280,280){$b$}
\put(820,1035){$c$}
\put(164,165){$\sigma$}  
\put(1164,152){$1$}       
\put(664,1022){$1$}       
\put(100,100){\arc{200}{0}{4.71}}     
\path(46,216)(100,200)(55,167)        
\spline(200,300)(277,733)(613,1020)   
\path(559,1007)(613,1020)(591,969)    
\spline(287,150)(700,0)(1113,150)     
\path(325,109)(287,150)(343,156)      
\put(1300,100){\arc{200}{4.71}{3.14}} 
\path(1345,167)(1300,200)(1354,216)     
\spline(650,983)(250,287)      
\path(297,318)(250,287)(253,343)      
\put(230,700){$_0$} 
\put(193,10){$_1$}  
\put(1155,10){$_0$}  
\put(680,77){$_1$}   
\put(455,585){$_{0,1}$}  
\end{picture}

\vspace{.3cm}

\hfill
\epsfig{file=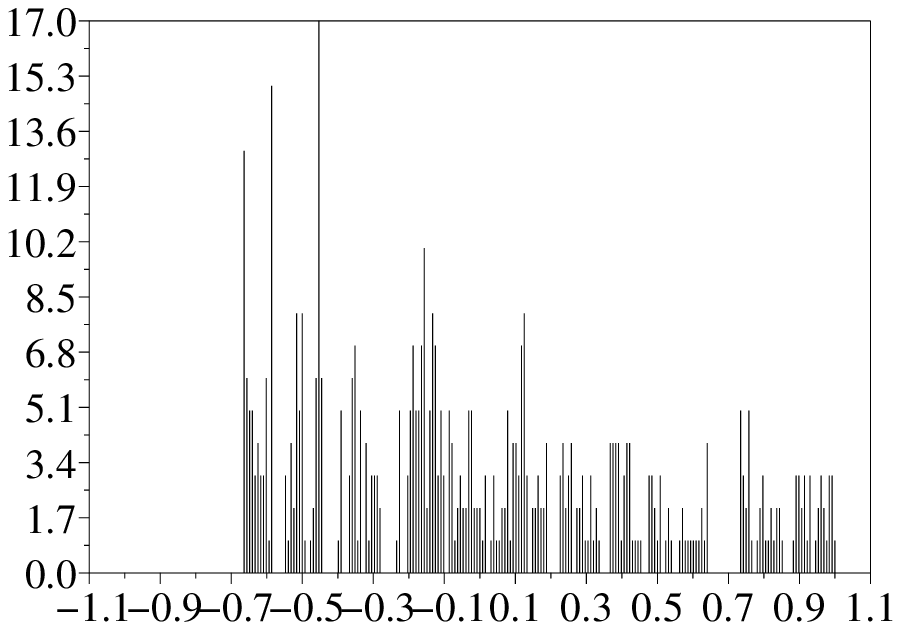,height=150pt}\\
\hline
\end{tabular}
\vspace{.8cm}

\noindent\begin{tabular}{@{}p{172pt}p{174pt}@{}}
\vspace{0.1cm}

\textbf{Automaton number $752$} \vspace{.2cm}

\begin{tabular}{@{}p{48pt}p{200pt}}

$a=\sigma(b,b)$

$b=(c,a)$

$c=(a,a)$& Group: Virtually $\Z^3$

Contracting: \textit{yes}

Self-replicating: \textit{no}\\
\end{tabular}

\begin{minipage}{230pt}
\setlength{\fboxrule}{0cm}
\begin{window}[5,r,{\fbox{\shortstack{\hspace{1.6cm}~\phantom{a}\\ \vspace{3.8cm}\\ \phantom{a}}}},{}]
Rels: $a^{2}$, $b^{2}$, $c^{2}$, $acbabacbab$, $acacbacacb$,\\
$abcabcacbacb$, $acbcbabacbcbab$, $abcbacbabcbacb$,\\
$acbcacbacbcacb$, $acacbcbacacbcb$, $abcbcabcacbcbacb$,\\
$acbcbcbabacbcbcbab$, $abcbcbacbabcbcbacb$.
\\
SF: $2^0$,$2^{1}$,$2^{3}$,$2^{5}$,$2^{7}$,$2^{8}$,$2^{10}$,$2^{11}$,$2^{13}$\\
Gr: 1,4,10,22,46,84,140,217,319,448\\
\end{window}
\end{minipage}
& \hfill~

\hfill
\begin{picture}(1450,1090)(0,130)
\put(200,200){\circle{200}} 
\put(1200,200){\circle{200}}
\put(700,1070){\circle{200}}
\allinethickness{0.2mm} \put(45,280){$a$} \put(1280,280){$b$}
\put(820,1035){$c$}
\put(164,165){$\sigma$}  
\put(1164,152){$1$}       
\put(664,1022){$1$}       
\put(300,200){\line(1,0){800}} 
\path(1050,225)(1100,200)(1050,175)   
\spline(287,150)(700,0)(1113,150)     
\path(325,109)(287,150)(343,156)      
\spline(750,983)(1150,287)     
\path(753,927)(750,983)(797,952)      
\spline(650,983)(250,287)      
\path(297,318)(250,287)(253,343)      
\put(650,250){$_{0,1}$} 
\put(890,585){$_0$} 
\put(680,77){$_1$}   
\put(455,585){$_{0,1}$}  
\end{picture}

\vspace{.3cm}

\hfill
\epsfig{file=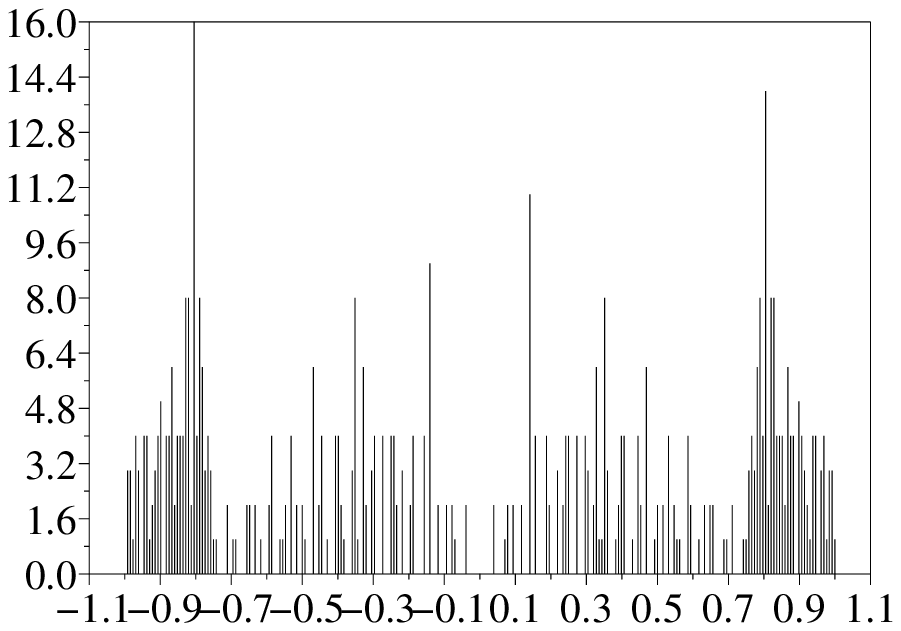,height=150pt}\\
\hline
\end{tabular}
\vspace{.8cm}

\noindent\begin{tabular}{@{}p{172pt}p{174pt}@{}} \textbf{Automata
number $775$ and $783$} \vspace{.2cm}

\begin{tabular}{@{}p{48pt}p{12pt}p{48pt}p{200pt}}

$a=\sigma(a,a)$

$b=(c,b)$

$c=(a,a)$&~

783:&

$a=\sigma(c,c)$

$b=(c,b)$

$c=(a,a)$&Group: $C_2\ltimes
IMG\left(\bigl(\frac{z-1}{z+1}\bigr)^2\right)$

Contracting: \emph{yes}

Self-replicating: \emph{yes}\\
\end{tabular}

\begin{minipage}{230pt}
\setlength{\fboxrule}{0cm}
\begin{window}[5,r,{\fbox{\shortstack{\hspace{1.6cm}~\phantom{a}\\ \vspace{3.8cm}\\ \phantom{a}}}},{}]
Rels: $a^{2}$, $b^{2}$, $c^{2}$, $acac$, $acbcbabcbcabcbabcb$,\\
$acbcbabcbacbcbabcb$, $abcbacbcbcacbcbcabcb$,\\
$acbcbacbcbabcbcabcbc$, $acbcbacbcbcacbcbcabcbc$
\\
SF: $2^0$,$2^{1}$,$2^{2}$,$2^{4}$,$2^{6}$,$2^{9}$,$2^{15}$,$2^{26}$,$2^{48}$\\
Gr: 1,4,9,17,30,51,85,140,229,367,579\\
Limit space:
\end{window}
\vspace{.1cm}
\noindent\epsfig{file=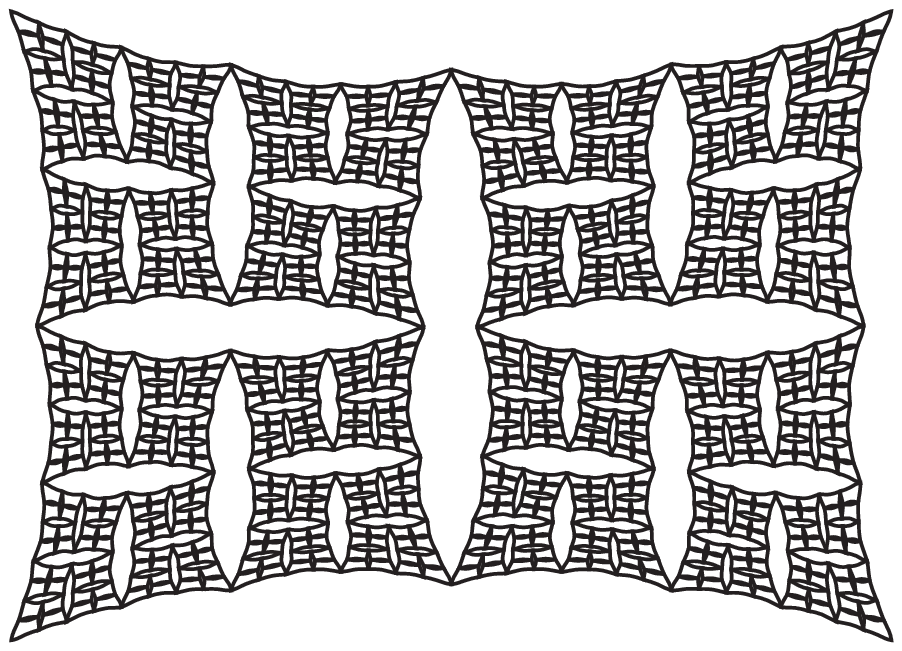,height=110pt}
\end{minipage} &
\hfill~

\hfill
\begin{picture}(1450,1090)(0,130)
\put(200,200){\circle{200}} 
\put(1200,200){\circle{200}}
\put(700,1070){\circle{200}}
\allinethickness{0.2mm} \put(45,280){$a$} \put(1280,280){$b$}
\put(820,1035){$c$}
\put(164,165){$\sigma$}  
\put(1164,152){$1$}       
\put(664,1022){$1$}       
\put(100,100){\arc{200}{0}{4.71}}     
\path(46,216)(100,200)(55,167)        
\put(1300,100){\arc{200}{4.71}{3.14}} 
\path(1345,167)(1300,200)(1354,216)     
\spline(750,983)(1150,287)     
\path(753,927)(750,983)(797,952)      
\spline(650,983)(250,287)      
\path(297,318)(250,287)(253,343)      
\put(190,10){$_{0,1}$}  
\put(890,585){$_0$} 
\put(1160,10){$_1$}  
\put(455,585){$_{0,1}$}  
\end{picture}

\vspace{.3cm}

\hfill\epsfig{file=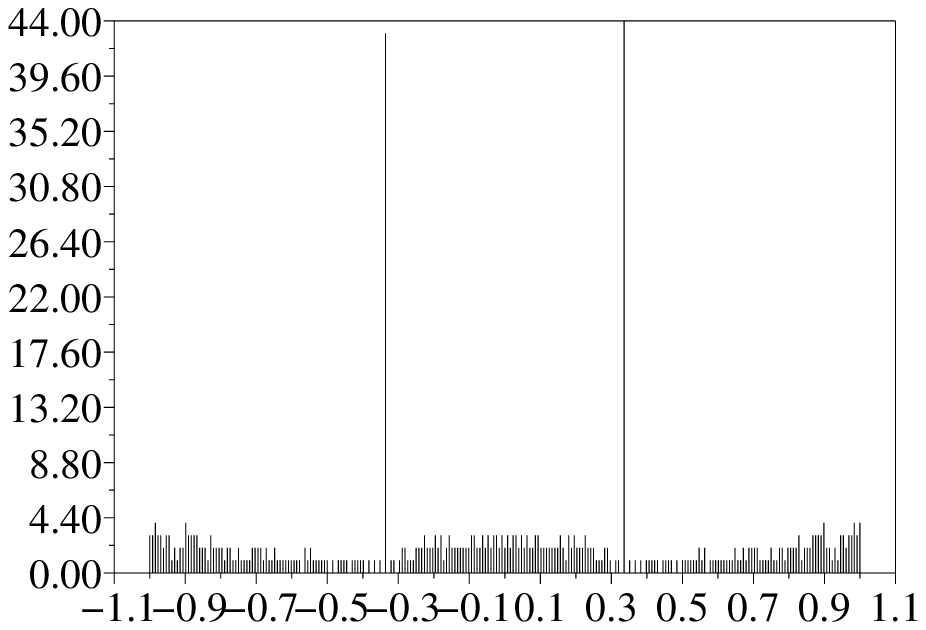,height=150pt}\\
\hline
\end{tabular}
\vspace{0.3cm}

\noindent\begin{tabular}{@{}p{172pt}p{174pt}@{}} \textbf{Automaton
number $802$} \vspace{.2cm}

\begin{tabular}{@{}p{48pt}p{200pt}}

$a=\sigma(a,a)$

$b=(c,c)$

$c=(a,a)$& Group: $C_2\times C_2\times C_2$

Contracting: \textit{yes}

Self-replicating: \textit{no}\\
\end{tabular}

\begin{minipage}{230pt}
\setlength{\fboxrule}{0cm}
\begin{window}[5,r,{\fbox{\shortstack{\hspace{1.6cm}~\phantom{a}\\ \vspace{3.8cm}\\ \phantom{a}}}},{}]
Rels: $a^{2}$, $b^{2}$, $c^{2}$, $abab$, $acac$, $bcbc$
\\
SF: $2^0$,$2^{1}$,$2^{2}$,$2^{3}$,$2^{3}$,$2^{3}$,$2^{3}$,$2^{3}$,$2^{3}$\\
Gr: 1,4,7,8,8,8,8,8,8,8,8\\
\end{window}
\end{minipage}
& \hfill~

\hfill
\begin{picture}(1450,1090)(0,130)
\put(200,200){\circle{200}} 
\put(1200,200){\circle{200}}
\put(700,1070){\circle{200}}
\allinethickness{0.2mm} \put(45,280){$a$} \put(1280,280){$b$}
\put(820,1035){$c$}
\put(164,165){$\sigma$}  
\put(1164,152){$1$}       
\put(664,1022){$1$}       
\put(100,100){\arc{200}{0}{4.71}}     
\path(46,216)(100,200)(55,167)        
\spline(750,983)(1150,287)     
\path(753,927)(750,983)(797,952)      
\spline(650,983)(250,287)      
\path(297,318)(250,287)(253,343)      
\put(190,10){$_{0,1}$}  
\put(820,585){$_{0,1}$} 
\put(455,585){$_{0,1}$}  
\end{picture}

\vspace{.3cm}

\hfill\epsfig{file=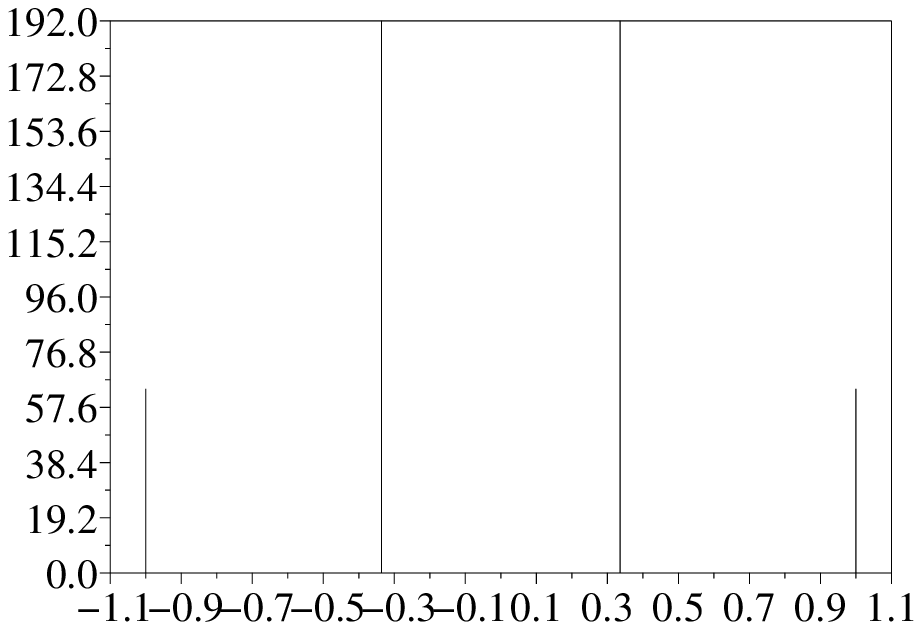,height=150pt}\\
\hline
\end{tabular}
\vspace{.3cm}

\noindent\begin{tabular}{@{}p{172pt}p{174pt}@{}} \textbf{Automaton
number $843$} \vspace{.2cm}

\begin{tabular}{@{}p{48pt}p{200pt}}

$a=\sigma(c,b)$

$b=(a,b)$

$c=(b,a)$& Group:

Contracting: \textit{no}

Self-replicating: \textit{yes}\\
\end{tabular}

\begin{minipage}{230pt}
\setlength{\fboxrule}{0cm}
\begin{window}[5,r,{\fbox{\shortstack{\hspace{1.6cm}~\phantom{a}\\ \vspace{3.8cm}\\ \phantom{a}}}},{}]
Rels: $acab^{-1}a^{-2}cab^{-1}aba^{-1}c^{-1}ba^{-1}c^{-1}$,\\
$acab^{-1}a^{-2}cb^{-1}ab^{-1}caba^{-1}c^{-2}ba^{-1}bc^{-1}$,\\
$acb^{-1}ab^{-1}ca^{-2}cab^{-1}ac^{-1}ba^{-1}bc^{-1}ba^{-1}c^{-1}$,\\
$acb^{-1}ab^{-1}ca^{-2}cb^{-1}ab^{-1}cac^{-1}ba^{-1}bc^{-2}ba^{-1}bc^{-1}$
\\
SF: $2^0$,$2^{1}$,$2^{3}$,$2^{5}$,$2^{8}$,$2^{14}$,$2^{24}$,$2^{43}$,$2^{81}$\\
Gr: 1,7,37,187,937,4687\\
\end{window}
\end{minipage}
& \hfill~

\hfill
\begin{picture}(1450,1090)(0,130)
\put(200,200){\circle{200}} 
\put(1200,200){\circle{200}}
\put(700,1070){\circle{200}}
\allinethickness{0.2mm} \put(45,280){$a$} \put(1280,280){$b$}
\put(820,1035){$c$}
\put(164,165){$\sigma$}  
\put(1164,152){$1$}       
\put(664,1022){$1$}       
\put(300,200){\line(1,0){800}} 
\path(1050,225)(1100,200)(1050,175)   
\spline(200,300)(277,733)(613,1020)   
\path(559,1007)(613,1020)(591,969)    
\spline(287,150)(700,0)(1113,150)     
\path(325,109)(287,150)(343,156)      
\put(1300,100){\arc{200}{4.71}{3.14}} 
\path(1345,167)(1300,200)(1354,216)     
\spline(650,983)(250,287)      
\path(297,318)(250,287)(253,343)      
\spline(1200,300)(1123,733)(787,1020) 
\path(1216,354)(1200,300)(1167,345)   
\put(230,700){$_0$} 
\put(680,240){$_1$} 
\put(680,77){$_0$}   
\put(1160,10){$_1$}  
\put(1115,700){$_0$}
\put(460,585){$_1$}  
\end{picture}

\vspace{.3cm}

\hfill
\epsfig{file=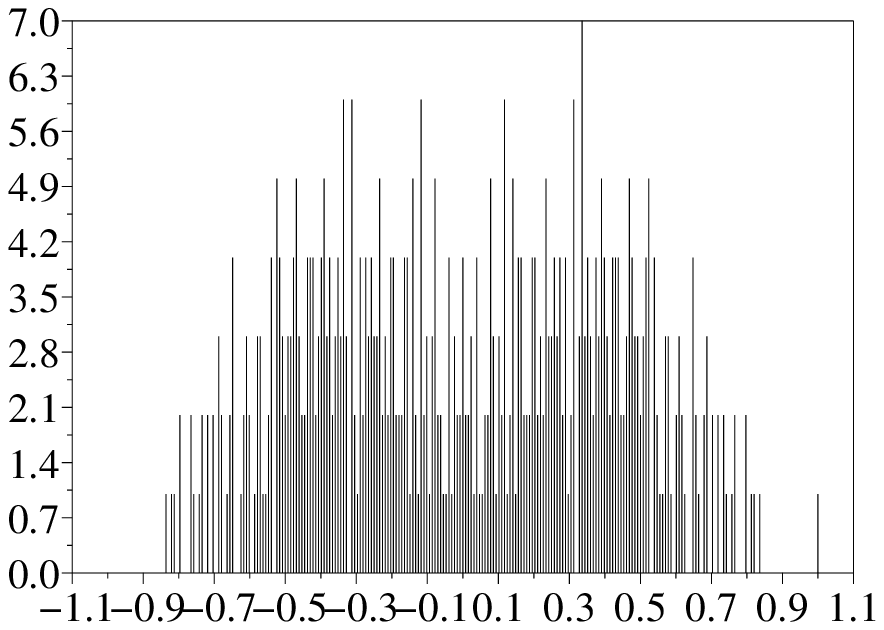,height=150pt}\\
\hline
\end{tabular}
\vspace{.8cm}

\noindent\begin{tabular}{@{}p{172pt}p{174pt}@{}} \textbf{Automaton
number $846$} \vspace{.2cm}

\begin{tabular}{@{}p{48pt}p{200pt}}

$a=\sigma(c,c)$

$b=(a,b)$

$c=(b,a)$& Group: $C_2\ast C_2\ast C_2$

Contracting: \textit{no}

Self-replicating: \textit{no}\\
\end{tabular}

\begin{minipage}{230pt}
\setlength{\fboxrule}{0cm}
\begin{window}[5,r,{\fbox{\shortstack{\hspace{1.6cm}~\phantom{a}\\ \vspace{3.8cm}\\ \phantom{a}}}},{}]
Rels: $a^{2}$, $b^{2}$, $c^{2}$
\\
SF: $2^0$,$2^{1}$,$2^{3}$,$2^{5}$,$2^{7}$,$2^{10}$,$2^{13}$,$2^{16}$,$2^{19}$\\
Gr: 1,4,10,22,46,94,190,382,766,1534\\
\end{window}
\end{minipage}
& \hfill~

\hfill
\begin{picture}(1450,1090)(0,130)
\put(200,200){\circle{200}} 
\put(1200,200){\circle{200}}
\put(700,1070){\circle{200}}
\allinethickness{0.2mm} \put(45,280){$a$} \put(1280,280){$b$}
\put(820,1035){$c$}
\put(164,165){$\sigma$}  
\put(1164,152){$1$}       
\put(664,1022){$1$}       
\spline(200,300)(277,733)(613,1020)   
\path(559,1007)(613,1020)(591,969)    
\spline(287,150)(700,0)(1113,150)     
\path(325,109)(287,150)(343,156)      
\put(1300,100){\arc{200}{4.71}{3.14}} 
\path(1345,167)(1300,200)(1354,216)     
\spline(650,983)(250,287)      
\path(297,318)(250,287)(253,343)      
\spline(1200,300)(1123,733)(787,1020) 
\path(1216,354)(1200,300)(1167,345)   
\put(150,700){$_{0,1}$} 
\put(680,77){$_0$}   
\put(1160,10){$_1$}  
\put(1115,700){$_0$}
\put(460,585){$_1$}  
\end{picture}

\vspace{.3cm}

\hfill\epsfig{file=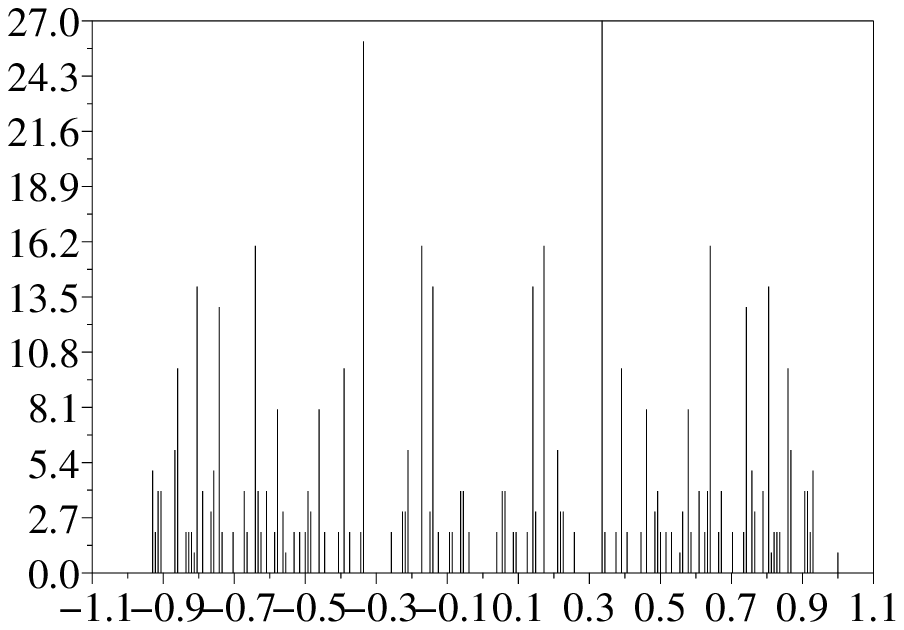,height=150pt}\\
\hline
\end{tabular}
\vspace{.3cm}

\noindent\begin{tabular}{@{}p{172pt}p{174pt}@{}} \textbf{Automaton
number $849$} \vspace{.2cm}

\begin{tabular}{@{}p{48pt}p{200pt}}

$a=\sigma(c,a)$

$b=(b,b)$

$c=(b,a)$& Group:

Contracting: \textit{no}

Self-replicating: \textit{yes}\\
\end{tabular}

\begin{minipage}{230pt}
\setlength{\fboxrule}{0cm}
\begin{window}[5,r,{\fbox{\shortstack{\hspace{1.6cm}~\phantom{a}\\ \vspace{3.8cm}\\ \phantom{a}}}},{}]
Rels: $b$, $a^{-1}c^{-1}ac^{-1}a^{-1}cac$,
$a^{-1}c^{-2}ac^{-1}a^{-1}c^{2}ac$,
$a^{-1}c^{-1}ac^{-2}a^{-1}cac^{2}$, $a^{-4}ca^{2}c^{-2}a^{2}c$,
$a^{-1}c^{-3}ac^{-1}a^{-1}c^{3}ac$,
$a^{-1}c^{-2}ac^{-2}a^{-1}c^{2}ac^{2}$,
$a^{-1}c^{-1}ac^{-3}a^{-1}cac^{3}$
\\
SF: $2^0$,$2^{1}$,$2^{3}$,$2^{6}$,$2^{12}$,$2^{23}$,$2^{45}$,$2^{88}$,$2^{174}$\\
Gr: 1,5,17,53,153,421,1125,2945,7589\\
\end{window}
\end{minipage}
& \hfill~

\hfill
\begin{picture}(1450,1090)(0,130)
\put(200,200){\circle{200}} 
\put(1200,200){\circle{200}}
\put(700,1070){\circle{200}}
\allinethickness{0.2mm} \put(45,280){$a$} \put(1280,280){$b$}
\put(820,1035){$c$}
\put(164,165){$\sigma$}  
\put(1164,152){$1$}       
\put(664,1022){$1$}       
\put(100,100){\arc{200}{0}{4.71}}     
\path(46,216)(100,200)(55,167)        
\spline(200,300)(277,733)(613,1020)   
\path(559,1007)(613,1020)(591,969)    
\put(1300,100){\arc{200}{4.71}{3.14}} 
\path(1345,167)(1300,200)(1354,216)     
\spline(650,983)(250,287)      
\path(297,318)(250,287)(253,343)      
\spline(1200,300)(1123,733)(787,1020) 
\path(1216,354)(1200,300)(1167,345)   
\put(230,700){$_0$} 
\put(193,10){$_1$}  
\put(1080,10){$_{0,1}$}  
\put(1115,700){$_0$}
\put(460,585){$_1$}  
\end{picture}

\vspace{.3cm}

\hfill\epsfig{file=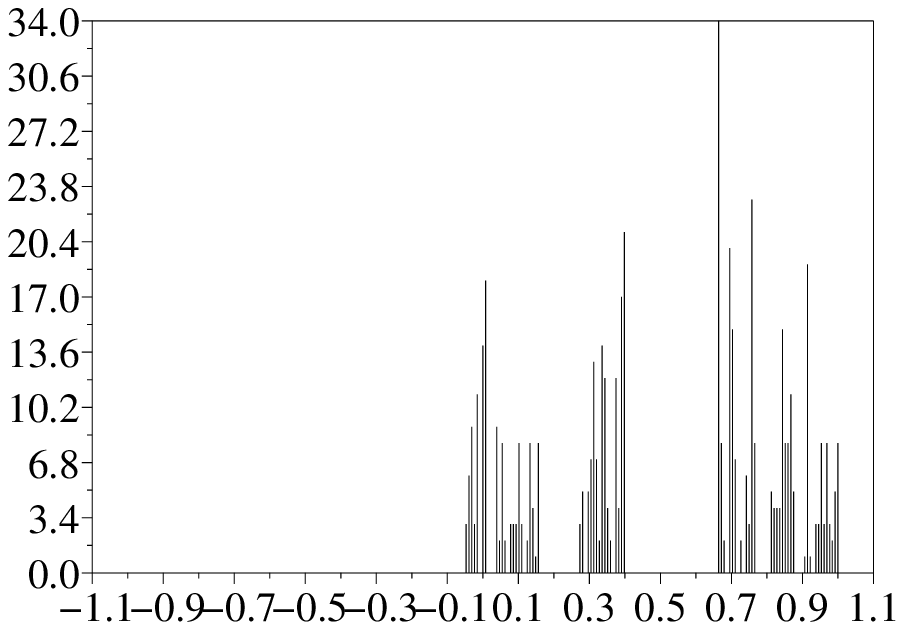,height=150pt}\\
\hline
\end{tabular}
\vspace{.3cm}

\noindent\begin{tabular}{@{}p{172pt}p{174pt}@{}}
\vspace{0.1cm}

\textbf{Automaton number $875$} \vspace{.2cm}

\begin{tabular}{@{}p{48pt}p{200pt}}

$a=\sigma(b,a)$

$b=(b,c)$

$c=(b,a)$& Group:

Contracting: \textit{no}

Self-replicating: \textit{yes}\\
\end{tabular}

\begin{minipage}{230pt}
\setlength{\fboxrule}{0cm}
\begin{window}[5,r,{\fbox{\shortstack{\hspace{1.6cm}~\phantom{a}\\ \vspace{3.8cm}\\ \phantom{a}}}},{}]
Rels: $a^{-1}ca^{-1}c$, $b^{-1}cb^{-1}c$,
$a^{-1}ba^{-1}ba^{-1}ba^{-1}b$,\\
$a^{-1}ba^{-1}bc^{-1}ac^{-1}ba^{-1}b$
\\
SF: $2^0$,$2^{1}$,$2^{3}$,$2^{7}$,$2^{13}$,$2^{25}$,$2^{47}$,$2^{90}$,$2^{176}$\\
Gr: 1,7,33,143,607,2563\\
\end{window}
\end{minipage}
& \hfill~

\hfill
\begin{picture}(1450,1090)(0,130)
\put(200,200){\circle{200}} 
\put(1200,200){\circle{200}}
\put(700,1070){\circle{200}}
\allinethickness{0.2mm} \put(45,280){$a$} \put(1280,280){$b$}
\put(820,1035){$c$}
\put(164,165){$\sigma$}  
\put(1164,152){$1$}       
\put(664,1022){$1$}       
\put(100,100){\arc{200}{0}{4.71}}     
\path(46,216)(100,200)(55,167)        
\put(300,200){\line(1,0){800}} 
\path(1050,225)(1100,200)(1050,175)   
\put(1300,100){\arc{200}{4.71}{3.14}} 
\path(1345,167)(1300,200)(1354,216)     
\spline(750,983)(1150,287)     
\path(753,927)(750,983)(797,952)      
\spline(650,983)(250,287)      
\path(297,318)(250,287)(253,343)      
\spline(1200,300)(1123,733)(787,1020) 
\path(1216,354)(1200,300)(1167,345)   
\put(680,240){$_0$} 
\put(193,10){$_1$}  
\put(1155,10){$_0$}  
\put(890,585){$_1$} 
\put(1115,700){$_0$}
\put(460,585){$_1$}  
\end{picture}

\vspace{.3cm}

\hfill
\epsfig{file=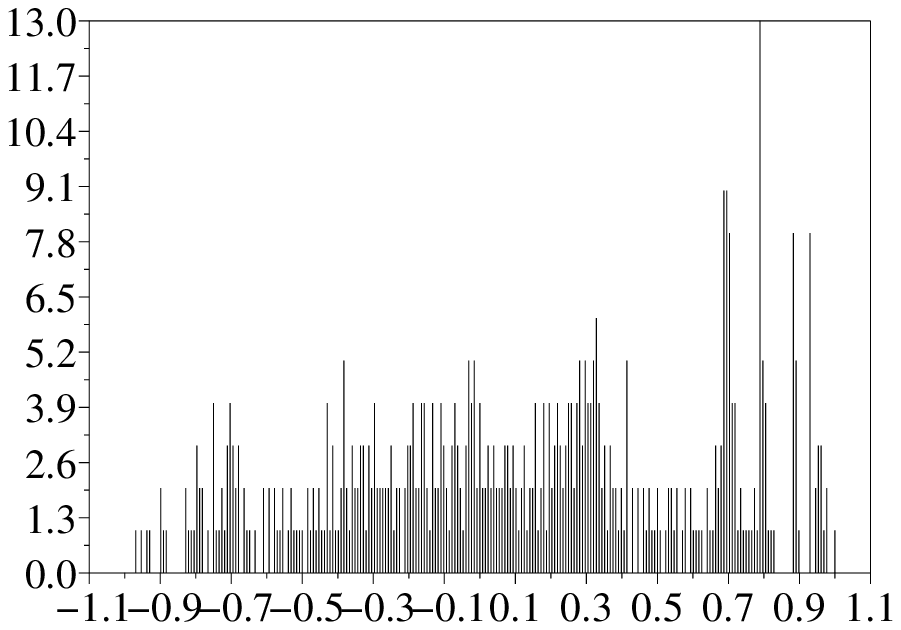,height=150pt}\\
\hline
\end{tabular}
\vspace{.8cm}

\noindent\begin{tabular}{@{}p{172pt}p{174pt}@{}} \textbf{Automaton
number $891$} \vspace{.2cm}

\begin{tabular}{@{}p{48pt}p{200pt}}

$a=\sigma(c,c)$

$b=(c,c)$

$c=(b,a)$& Group: $C_2\ltimes Lampighter$

Contracting: \textit{no}

Self-replicating: \textit{yes}\\
\end{tabular}

\begin{minipage}{230pt}
\setlength{\fboxrule}{0cm}
\begin{window}[5,r,{\fbox{\shortstack{\hspace{1.6cm}~\phantom{a}\\ \vspace{3.8cm}\\ \phantom{a}}}},{}]
Rels: $a^{2}$, $b^{2}$, $c^{2}$, $abab$, $acabcbacbacb$,\\
$acabcabcbacb$, $acacabcacbacacbacacb$,\\ $acacabcacabcacbacacb$,
$acacabcbcbacbcbacacb$,\\ $acacabcbcabcbcbacacb$,
$acabcbcabcacbacbcbac$
\\
SF: $2^0$,$2^{1}$,$2^{3}$,$2^{6}$,$2^{7}$,$2^{9}$,$2^{10}$,$2^{11}$,$2^{12}$\\
Gr: 1,4,9,17,30,51,82,128,198,304,456\\
\end{window}
\end{minipage}
& \hfill~

\hfill
\begin{picture}(1450,1090)(0,130)
\put(200,200){\circle{200}} 
\put(1200,200){\circle{200}}
\put(700,1070){\circle{200}}
\allinethickness{0.2mm} \put(45,280){$a$} \put(1280,280){$b$}
\put(820,1035){$c$}
\put(164,165){$\sigma$}  
\put(1164,152){$1$}       
\put(664,1022){$1$}       
\spline(200,300)(277,733)(613,1020)   
\path(559,1007)(613,1020)(591,969)    
\spline(750,983)(1150,287)     
\path(753,927)(750,983)(797,952)      
\spline(650,983)(250,287)      
\path(297,318)(250,287)(253,343)      
\spline(1200,300)(1123,733)(787,1020) 
\path(1216,354)(1200,300)(1167,345)   
\put(150,700){$_{0,1}$} 
\put(820,585){$_{0,1}$} 
\put(1115,700){$_0$}
\put(460,585){$_1$}  
\end{picture}

\vspace{.3cm}

\hfill\epsfig{file=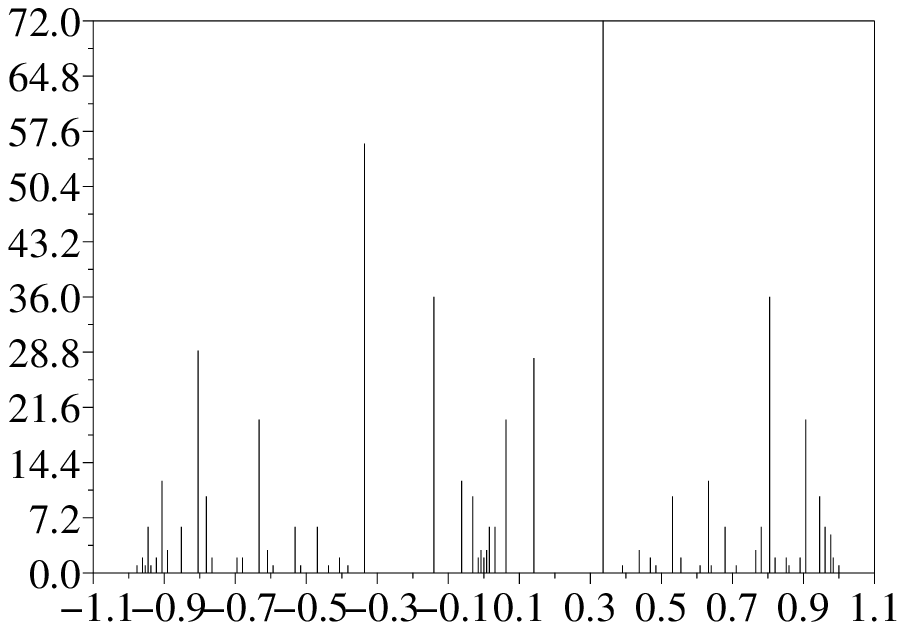,height=150pt}\\
\hline
\end{tabular}
\vspace{.3cm}

\noindent\begin{tabular}{@{}p{172pt}p{174pt}@{}}
\vspace{0.1cm}

\textbf{Automaton number $2193$} \vspace{.2cm}

\begin{tabular}{@{}p{48pt}p{200pt}}

$a=\sigma(c,b)$

$b=\sigma(a,a)$

$c=(a,a)$& Group: contains $Lamplighter$ group

Contracting: \textit{no}

Self-replicating: \textit{yes}\\
\end{tabular}

\begin{minipage}{230pt}
\setlength{\fboxrule}{0cm}
\begin{window}[5,r,{\fbox{\shortstack{\hspace{1.6cm}~\phantom{a}\\ \vspace{3.8cm}\\ \phantom{a}}}},{}]
Rels: $b^{-1}c^{-1}bc$, $b^{-2}c^{2}$, $b^{-1}cb^{-1}c$, $a^{-4}$,\\
$b^{-4}$, $b^{-1}c^{-2}b^{-1}$, $b^{-1}c^{-1}b^{-1}c^{-1}$,
$a^{-2}ba^{2}b$,\\
$a^{-1}b^{-1}c^{-1}a^{-1}cb$,
 $a^{-2}ba^{-2}b$, $a^{-2}ca^{2}c$,
$a^{-2}ca^{-2}c$, $b^{-1}c^{-1}a^{-1}bca^{-1}$,
$a^{-1}b^{-1}a^{-2}b^{-1}a^{-1}$.
\\
SF: $2^0$,$2^{1}$,$2^{3}$,$2^{6}$,$2^{7}$,$2^{9}$,$2^{10}$,$2^{11}$,$2^{12}$\\
Gr: 1,7,27,65,120,204,328,512,792,1216\\
\end{window}
\end{minipage}
& \hfill~

\hfill
\begin{picture}(1450,1090)(0,130)
\put(200,200){\circle{200}} 
\put(1200,200){\circle{200}}
\put(700,1070){\circle{200}}
\allinethickness{0.2mm} \put(45,280){$a$} \put(1280,280){$b$}
\put(820,1035){$c$}
\put(164,165){$\sigma$}  
\put(1164,165){$\sigma$}  
\put(664,1022){$1$}       
\put(300,200){\line(1,0){800}} 
\path(1050,225)(1100,200)(1050,175)   
\spline(200,300)(277,733)(613,1020)   
\path(559,1007)(613,1020)(591,969)    
\spline(287,150)(700,0)(1113,150)     
\path(325,109)(287,150)(343,156)      
\spline(650,983)(250,287)      
\path(297,318)(250,287)(253,343)      
\put(230,700){$_0$} 
\put(680,240){$_1$} 
\put(650,87){$_{0,1}$}   
\put(455,585){$_{0,1}$}  
\end{picture}

\vspace{.3cm}

\hfill
\epsfig{file=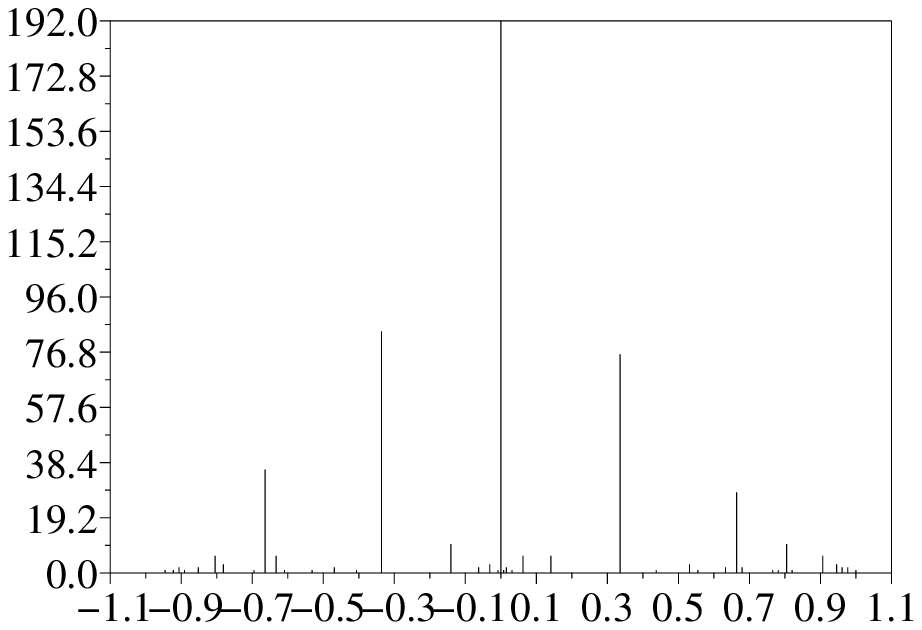,height=150pt}\\
\hline
\end{tabular}
\vspace{.8cm}

\noindent\begin{tabular}{@{}p{172pt}p{174pt}@{}}
\vspace{0.1cm}

\textbf{Automaton number $2280$} \vspace{.2cm}

\begin{tabular}{@{}p{48pt}p{200pt}}

$a=\sigma(c,a)$

$b=\sigma(b,a)$

$c=(b,a)$& Group:

Contracting: \textit{no}

Self-replicating: \textit{yes}\\
\end{tabular}

\begin{minipage}{230pt}
\setlength{\fboxrule}{0cm}
\begin{window}[5,r,{\fbox{\shortstack{\hspace{1.6cm}~\phantom{a}\\ \vspace{3.8cm}\\ \phantom{a}}}},{}]
Rels: $a^{-1}ba^{-1}b$, $b^{-1}cb^{-1}c$, $b^{-1}ca^{-1}ba^{-1}c$,\\
$a^{-2}b^{2}a^{-1}b^{-1}ab$, $a^{-2}bab^{-2}ab$,
$b^{-2}c^{2}b^{-1}c^{-1}bc$,\\ $b^{-2}cbc^{-2}bc$,
$a^{-1}ca^{-1}ca^{-1}ca^{-1}c$, $b^{-1}ca^{-1}cb^{-1}ca^{-1}c$,\\
$a^{-1}bab^{-1}a^{-1}b^{2}a^{-1}$,
$b^{-1}cbc^{-1}b^{-1}c^{2}b^{-1}$,\\
$a^{-2}bac^{-1}bc^{-1}b^{-1}ab$.
\\
SF: $2^0$,$2^{1}$,$2^{3}$,$2^{7}$,$2^{13}$,$2^{25}$,$2^{47}$,$2^{90}$,$2^{176}$\\
Gr: 1,7,33,143,597,2465\\
\end{window}
\end{minipage}
& \hfill~

\hfill
\begin{picture}(1450,1090)(0,130)
\put(200,200){\circle{200}} 
\put(1200,200){\circle{200}}
\put(700,1070){\circle{200}}
\allinethickness{0.2mm} \put(45,280){$a$} \put(1280,280){$b$}
\put(820,1035){$c$}
\put(164,165){$\sigma$}  
\put(1164,165){$\sigma$}  
\put(664,1022){$1$}       
\put(100,100){\arc{200}{0}{4.71}}     
\path(46,216)(100,200)(55,167)        
\spline(200,300)(277,733)(613,1020)   
\path(559,1007)(613,1020)(591,969)    
\spline(287,150)(700,0)(1113,150)     
\path(325,109)(287,150)(343,156)      
\put(1300,100){\arc{200}{4.71}{3.14}} 
\path(1345,167)(1300,200)(1354,216)     
\spline(650,983)(250,287)      
\path(297,318)(250,287)(253,343)      
\spline(1200,300)(1123,733)(787,1020) 
\path(1216,354)(1200,300)(1167,345)   
\put(230,700){$_0$} 
\put(193,10){$_1$}  
\put(1155,10){$_0$}  
\put(680,77){$_1$}   
\put(1115,700){$_0$}
\put(460,585){$_1$}  
\end{picture}

\vspace{.3cm}

\hfill
\epsfig{file=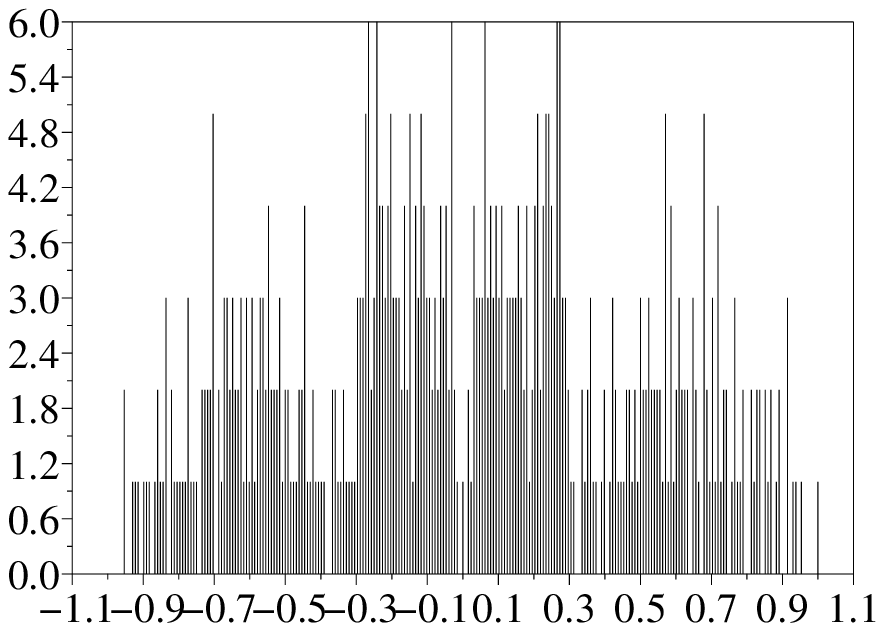,height=150pt}\\
\hline
\end{tabular}
\vspace{.8cm}

\noindent\begin{tabular}{@{}p{172pt}p{174pt}@{}} \textbf{Automaton
number $2294$} \vspace{.2cm}

\begin{tabular}{@{}p{48pt}p{200pt}}

$a=\sigma(b,c)$

$b=\sigma(c,a)$

$c=(b,a)$& Group: \textit{Baumslag-Solitar group $B(1,-3)$}

Contracting: \textit{no}

Self-replicating: \textit{yes}\\
\end{tabular}

\begin{minipage}{230pt}
\setlength{\fboxrule}{0cm}
\begin{window}[5,r,{\fbox{\shortstack{\hspace{1.6cm}~\phantom{a}\\ \vspace{3.8cm}\\ \phantom{a}}}},{}]
Rels: $b^{-1}ca^{-1}c$, $(ca^{-1})^a(ca^{-1})^3$
\\
SF: $2^0$,$2^{1}$,$2^{2}$,$2^{4}$,$2^{6}$,$2^{8}$,$2^{10}$,$2^{12}$,$2^{14}$\\
Gr: 1,7,33,127,433,1415\\
\end{window}
\end{minipage}
& \hfill~

\hfill
\begin{picture}(1450,1090)(0,130)
\put(200,200){\circle{200}} 
\put(1200,200){\circle{200}}
\put(700,1070){\circle{200}}
\allinethickness{0.2mm} \put(45,280){$a$} \put(1280,280){$b$}
\put(820,1035){$c$}
\put(164,165){$\sigma$}  
\put(1164,165){$\sigma$}  
\put(664,1022){$1$}       
\put(300,200){\line(1,0){800}} 
\path(1050,225)(1100,200)(1050,175)   
\spline(200,300)(277,733)(613,1020)   
\path(559,1007)(613,1020)(591,969)    
\spline(287,150)(700,0)(1113,150)     
\path(325,109)(287,150)(343,156)      
\spline(750,983)(1150,287)     
\path(753,927)(750,983)(797,952)      
\spline(650,983)(250,287)      
\path(297,318)(250,287)(253,343)      
\spline(1200,300)(1123,733)(787,1020) 
\path(1216,354)(1200,300)(1167,345)   
\put(680,240){$_0$} 
\put(230,700){$_1$} 
\put(890,585){$_0$} 
\put(680,77){$_1$}   
\put(1115,700){$_0$}
\put(460,585){$_1$}  
\end{picture}

\vspace{.3cm}

\hfill\epsfig{file=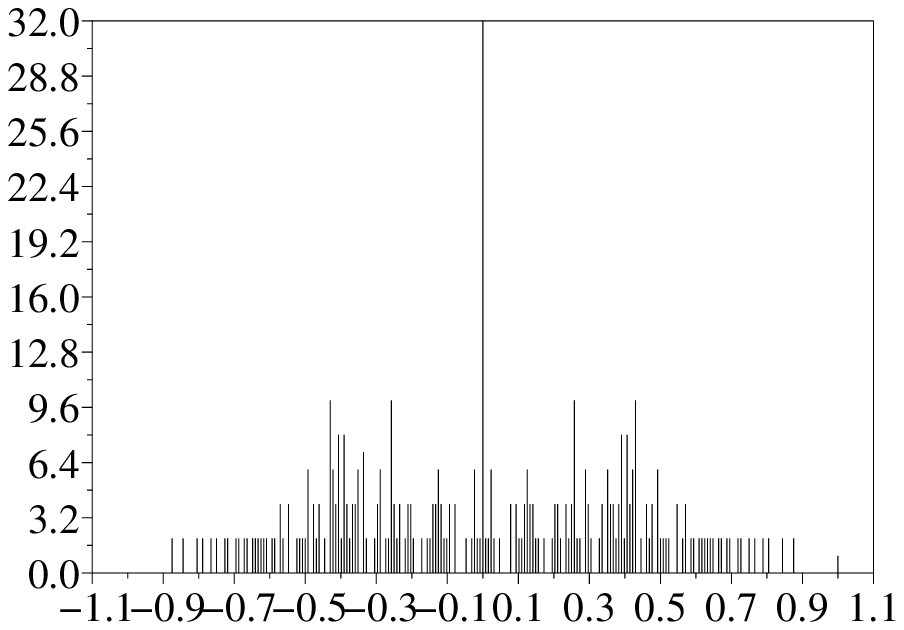,height=150pt}\\
\hline
\end{tabular}

\noindent\begin{tabular}{@{}p{172pt}p{174pt}@{}}
\vspace{0.1cm}

\textbf{Automaton number $2396$} \vspace{.2cm}

\begin{tabular}{@{}p{48pt}p{200pt}}

$a=\sigma(b,a)$

$b=\sigma(c,b)$

$c=(c,a)$& Group: \textit{Boltenkov group}

Contracting: \textit{no}

Self-replicating: \textit{yes}\\
\end{tabular}

\begin{minipage}{230pt}
\setlength{\fboxrule}{0cm}
\begin{window}[5,r,{\fbox{\shortstack{\hspace{1.6cm}~\phantom{a}\\ \vspace{3.8cm}\\ \phantom{a}}}},{}]
Rels: $acb^{-1}ca^{-2}cb^{-1}cac^{-1}bc^{-2}bc^{-1}$,\\
$acb^{-1}ca^{-2}cb^{-1}a^2c^{-1}b^{-1}a^2c^{-1}bc^{-1}a^{-1}bca^{-2}bc^{-1}$,\\
$acb^{-1}a^2c^{-1}b^{-1}a^{-1}cb^{-1}cbca^{-2}bc^{-2}bc^{-1}$,\\
$bcb^{-1}ca^{-1}b^{-1}cb^{-1}a^2c^{-1}ac^{-1}ba^{-2}bc^{-1}$
\\
SF: $2^0$,$2^{1}$,$2^{3}$,$2^{6}$,$2^{12}$,$2^{24}$,$2^{46}$,$2^{90}$,$2^{176}$\\
Gr: 1,7,37,187,937,4687\\
\end{window}
\end{minipage}
& \hfill~

\hfill
\begin{picture}(1450,1090)(0,130)
\put(200,200){\circle{200}} 
\put(1200,200){\circle{200}}
\put(700,1070){\circle{200}}
\allinethickness{0.2mm} \put(45,280){$a$} \put(1280,280){$b$}
\put(820,1035){$c$}
\put(164,165){$\sigma$}  
\put(1164,165){$\sigma$}  
\put(664,1022){$1$}       
\put(100,100){\arc{200}{0}{4.71}}     
\path(46,216)(100,200)(55,167)        
\put(300,200){\line(1,0){800}} 
\path(1050,225)(1100,200)(1050,175)   
\put(1300,100){\arc{200}{4.71}{3.14}} 
\path(1345,167)(1300,200)(1354,216)     
\spline(750,983)(1150,287)     
\path(753,927)(750,983)(797,952)      
\spline(650,983)(250,287)      
\path(297,318)(250,287)(253,343)      
\put(700,1211){\arc{200}{2.36}{0.79}} 
\path(820,1168)(771,1141)(779,1196)   
\put(680,240){$_0$} 
\put(193,10){$_1$}  
\put(890,585){$_0$} 
\put(1160,10){$_1$}  
\put(545,1261){$_0$}  
\put(460,585){$_1$}  
\end{picture}

\vspace{.3cm}

\hfill
\epsfig{file=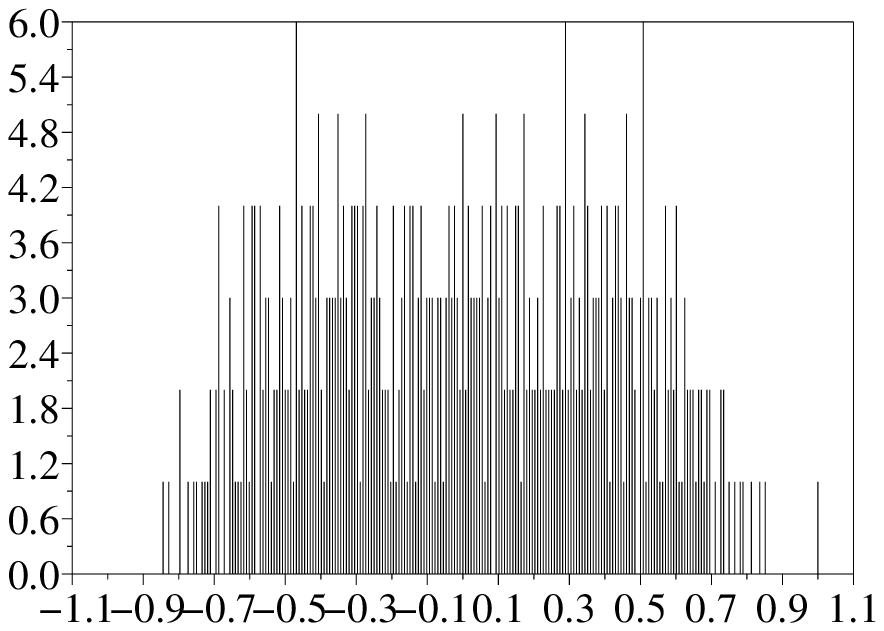,height=150pt}\\
\hline
\end{tabular}
\vspace{.8cm}

\noindent\begin{tabular}{@{}p{172pt}p{174pt}@{}}
\vspace{0.1cm}

\textbf{Automaton number $2398$} \vspace{.2cm}

\begin{tabular}{@{}p{48pt}p{200pt}}

$a=\sigma(a,b)$

$b=\sigma(c,b)$

$c=(c,a)$& Group: \textit{Dahmani Group}

Contracting: \textit{no}

Self-replicating: \textit{yes}\\
\end{tabular}

\begin{minipage}{230pt}
\setlength{\fboxrule}{0cm}
\begin{window}[5,r,{\fbox{\shortstack{\hspace{1.6cm}~\phantom{a}\\ \vspace{3.8cm}\\ \phantom{a}}}},{}]
Rels: $cba$, $b^{-1}a^{-1}b^{2}ca^{2}$, $a^{-2}c^{-1}acb^{-1}ab$,\\
$a^{-1}c^{-1}b^{-1}acaba^{-1}$,
$b^{-1}a^{2}b^{-1}a^{-1}b^{2}a^{-1}$,\\
$c^{-1}b^{-3}c^{-1}ba^{-2}$,
$c^{-1}b^{-1}a^{-1}b^{2}cab^{-1}$,\\
$c^{-1}ab^{-1}a^{-1}b^{2}cb^{-1}$, $a^{-3}c^{-1}b^{-2}ac^{-1}$,\\
$ca^{-1}bc^{-1}a^{-1}ca^{3}$, $cabcba^{-1}c^{-1}b^{-1}a$,\\
$ab^{-3}c^{-1}ba^{-2}b$, $ba^{-1}c^{-1}acb^{-1}abc$,\\
$bc^{-1}b^{-1}acaba^{-1}c$
\\
SF: $2^0$,$2^{1}$,$2^{3}$,$2^{6}$,$2^{12}$,$2^{23}$,$2^{45}$,$2^{88}$,$2^{174}$\\
Gr: 1,7,31,127,483,1823\\
\end{window}
\end{minipage}
& \hfill~

\hfill
\begin{picture}(1450,1090)(0,130)
\put(200,200){\circle{200}} 
\put(1200,200){\circle{200}}
\put(700,1070){\circle{200}}
\allinethickness{0.2mm} \put(45,280){$a$} \put(1280,280){$b$}
\put(820,1035){$c$}
\put(164,165){$\sigma$}  
\put(1164,165){$\sigma$}  
\put(664,1022){$1$}       
\put(100,100){\arc{200}{0}{4.71}}     
\path(46,216)(100,200)(55,167)        
\put(300,200){\line(1,0){800}} 
\path(1050,225)(1100,200)(1050,175)   
\put(1300,100){\arc{200}{4.71}{3.14}} 
\path(1345,167)(1300,200)(1354,216)     
\spline(750,983)(1150,287)     
\path(753,927)(750,983)(797,952)      
\spline(650,983)(250,287)      
\path(297,318)(250,287)(253,343)      
\put(700,1211){\arc{200}{2.36}{0.79}} 
\path(820,1168)(771,1141)(779,1196)   
\put(190,10){$_0$}  
\put(680,240){$_1$} 
\put(890,585){$_0$} 
\put(1160,10){$_1$}  
\put(545,1261){$_0$}  
\put(460,585){$_1$}  
\end{picture}

\vspace{.3cm}

\hfill
\epsfig{file=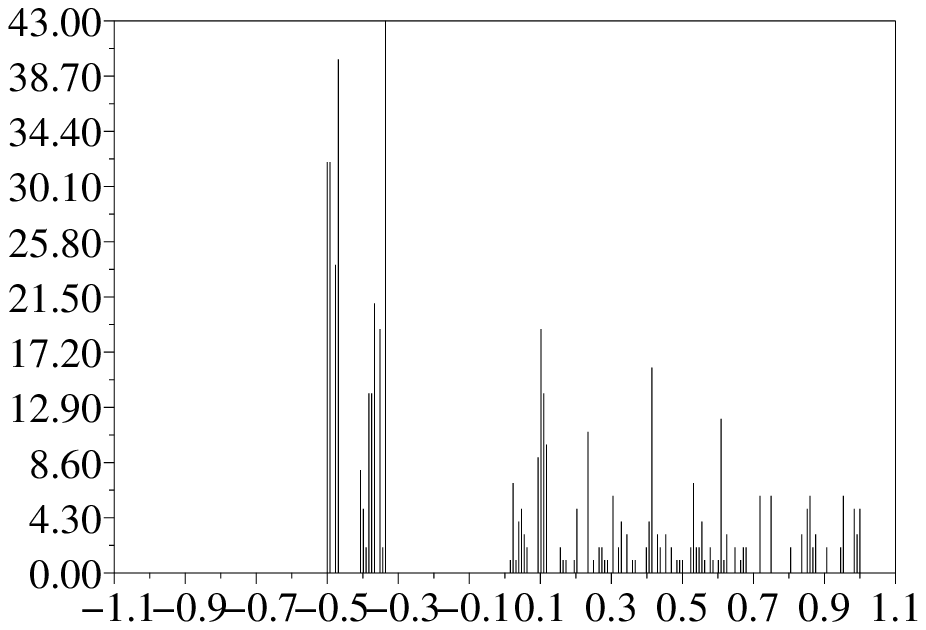,height=150pt}\\
\hline
\end{tabular}
\vspace{.8cm}

\noindent\begin{tabular}{@{}p{172pt}p{174pt}@{}} \textbf{Automaton
number $2852$} \vspace{.2cm}

\begin{tabular}{@{}p{48pt}p{200pt}}

$a=\sigma(b,c)$

$b=\sigma(b,a)$

$c=(c,c)$& Group: Isomorphic to $G_{849}$

Contracting: \textit{no}

Self-replicating: \textit{yes}\\
\end{tabular}

\begin{minipage}{230pt}
\setlength{\fboxrule}{0cm}
\begin{window}[5,r,{\fbox{\shortstack{\hspace{1.6cm}~\phantom{a}\\ \vspace{3.8cm}\\ \phantom{a}}}},{}]
Rels: $c$, $a^{-2}bab^{-2}ab$, $a^{-2}ba^{-1}bab^{-2}ab^{-1}ab$,
$a^{-1}bab^{-1}ab^{-2}aba^{-1}ba^{-1}$
\\
SF: $2^0$,$2^{1}$,$2^{3}$,$2^{6}$,$2^{12}$,$2^{23}$,$2^{45}$,$2^{88}$,$2^{174}$\\
Gr: 1,5,17,53,153,429,1189\\
\end{window}
\end{minipage}
& \hfill~

\hfill
\begin{picture}(1450,1090)(0,130)
\put(200,200){\circle{200}} 
\put(1200,200){\circle{200}}
\put(700,1070){\circle{200}}
\allinethickness{0.2mm} \put(45,280){$a$} \put(1280,280){$b$}
\put(820,1035){$c$}
\put(164,165){$\sigma$}  
\put(1164,165){$\sigma$}  
\put(664,1022){$1$}       
\put(300,200){\line(1,0){800}} 
\path(1050,225)(1100,200)(1050,175)   
\spline(200,300)(277,733)(613,1020)   
\path(559,1007)(613,1020)(591,969)    
\spline(287,150)(700,0)(1113,150)     
\path(325,109)(287,150)(343,156)      
\put(1300,100){\arc{200}{4.71}{3.14}} 
\path(1345,167)(1300,200)(1354,216)     
\put(700,1211){\arc{200}{2.36}{0.79}} 
\path(820,1168)(771,1141)(779,1196)   
\put(680,240){$_0$} 
\put(230,700){$_1$} 
\put(1155,10){$_0$}  
\put(680,77){$_1$}   
\put(465,1261){$_{0,1}$}  
\end{picture}

\vspace{.3cm}

\hfill\epsfig{file=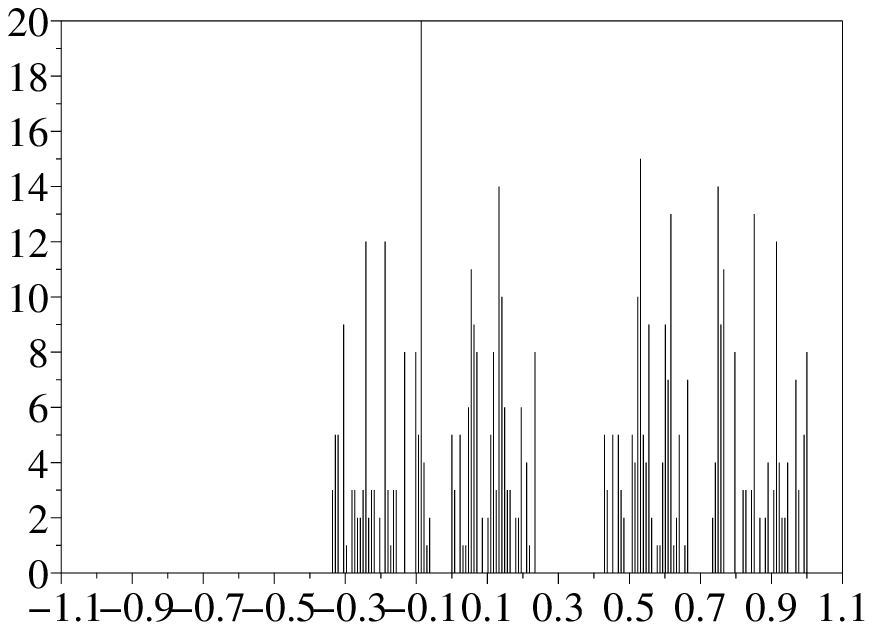,height=150pt}\\
\hline
\end{tabular}

\section{Proofs of some facts about the selected
groups}\label{proofs}

\noindent\textbf{741}: For $G_{741}$, we have $a=\sigma(c,a)$,
$b=(b,a)$ and $c=(a,a)$.

The states $a$ and $c$ form a $2$-state automaton generating $\Z$
(see Theorem~\ref{thm:class22} and its proof
in~\cite{gns00:automata}). Since $b^n = (b^n,a^n)$ we see that $b$
also has infinite order and that $G_{741}$ is not contracting ($b^n$
would have to belong to the nucleus for all $n$).

\medskip

\noindent\textbf{752}: For $G_{752}$, we have $a=\sigma(b,b)$,
$b=(c,a)$ and $c=(a,a)$.

We claim that this group contains $\Z^3$ as a subgroup of index $4$.
It is contracting with nucleus consisting of $41$ elements.

Since $a^2 = (b^2,b^2)$, $b^2 = (c^2,a^2)$ and $c^2 = (a^2,a^2)$,
all generators have order~2.

Let $x=ca$, $y=babc$, and consider the subgroup $K=\langle x,y
\rangle$. Direct computations show that $x$ and $y$ commute
($xy=((cbab)^{ca},abcb)=((y^{-1})^x,abcb)$ and
$yx=(cbab,abcb)=(y^{-1},abcb)$). Conjugating by
$\gamma=(\gamma,bc\gamma)$ leads to the self-similar copy $K'$ of
$K$ generated by $x' = \sigma((y')^{-1},(x')^{-1})$ and $y' =
\sigma((y')^{-1}x',1)$, where $x' = x^\gamma$ and $y'=y^\gamma$.
Since $(x')^2 = ((x')^{-1}(y')^{-1},(y')^{-1}(x')^{-1})$ and $(y')^2
= ((y')^{-1}x',(y')^{-1}x')$, the virtual endomorphism of $K'$ is
given by the matrix
\[
 A=\left(\begin{array}{cc}
               -\frac12&\frac12\\
               -\frac12&-\frac12\\
 \end{array} \right).
\]
The eigenvalues $\lambda=-\frac12\pm\frac12i$ of this matrix are not
algebraic integers, hence, according
to~\cite{nekrash_s:12endomorph}, the group $K'$ is free abelian of
rank $2$, and so is $K$.

Since all generators have order $2$, the subgroup $H=\langle
ba,cb\rangle$ has index $2$ in $G_{752}$. The stabilizer $\St_H(1)$
of the first level has index $2$ in $H$ and is generated by $cb$,
$cb^{ba}$ and $(ba)^2$. Conjugating these three generators by
$g=(1,b)$, we obtain
\begin{alignat*}{3}
g_1 &= \left(cb\right)^g          &&= (x^{-1},&&1),\\
g_2 &= \left((cb)^{ba}\right)^{g} &&=(1,&&x),\\
g_3 &= \left((ba)^2\right)^{g}    &&=(y^{-1},&&y).
\end{alignat*}
Therefore, $g_1$, $g_2$, and $g_3$ commute. If
$g_1^{n_1}g_2^{n_2}g_3^{n_3}=1$, then all sections must be trivial,
implying that $x^{-n_1}y^{-n_3}=x^{n_2}y^{n_3}=1$. But $K$ is free
abelian, whence $n_i=0$, $i=1,2,3$. Thus, $\St_H(1)$ is a free
abelian group of rank $3$.

\medskip

\noindent\textbf{775}: For $G_{775}$, we have $a=\sigma(a,a)$,
$b=(c,b)$ and $c=(a,a)$.

This group was mentioned in the previous
paper~\cite{bondarenko_gkmnss:clas32}, where it was proved that it
is conjugate to $G_{783}$ and isomorphic to $G_{2205}$. Some
additional properties of this group were mentioned without a proof.
Here we list and prove some of those properties.

Conjugating the generators by $g = \sigma(g,g)$ we obtain
\[
 a' = \s (a',a'), \qquad b' = (b',c'), \qquad c' = (a',a'),
\]
where $a'=a^g$, $b'=b^g$ and $c'=c^g$.  The tree automorphisms $a'$,
$b'$ and $c'$ are precisely the generators of the group $G_{793}$
and we choose to work with this group instead of $G_{775}$. In fact,
from now on the generators of $G_{793}$ will be denoted by $a$, $b$
and $c$ and the group $G_{793}$ will be denoted by $G$. Thus
\[
 a = \s (a,a), \qquad b = (b,c), \qquad c = (a,a)
\]
and $G = \langle a,b,c\rangle$.

Since $a^2 = (a^2,a^2)$, $b^2 = (b^2,c^2)$ and $c^2 = (a^2,a^2)$,
all three generators have order 2. Moreover, since $ac = \s$ has
order 2, $a$ and $c$ commute, i.e., $\langle a,c \rangle = C_2
\times C_2$ is the Klein Viergruppe. Denote this subgroup of $G$ by
$A$.

We note that the element $x=ba$ has infinite order. Indeed,
\[ x = ba = (b,c)\s(a,a) = \s(ca,ba) = \s(\s,x). \]
Therefore $x^2=(x \s, \s x) = (x,\s,\s,x )$. Since $x^2$ fixes $00$,
and has $x$ as a section, $x$ must have infinite order (the element
$x$ cannot have odd order since it acts nontrivially at the root; on
the other hand $x^{2n}=1$ implies $x^n=1$, so $x$ cannot have even
order as well).

\begin{prop}
The subgroup $H=\langle x,y \rangle$ of $G$, where $x=ba$ and
$y=cabc$ is torsion free. Moreover
\[ G = A \ltimes H = (C_2 \times C_2) \ltimes H \]
\end{prop}

\begin{proof}
The equalities
\begin{align*}
 & x^a = abaa = ab = x^{-1} \qquad & y^a = acabca = cbac = y^{-1} \\
 & x^b = bbab = ab = x^{-1} \qquad & y^b = bcabcb = bacbacab = xy^{-1}x^{-1}\\
 & x^c = cbac = y^{-1}      \qquad & y^c = ccabcc = ab = x^{-1}
\end{align*}
show that $H$ is normal in $G$ (in fact they even show that $H$ is
the normal closure of $x$ in $G$). Since $\{x,y,a,c\}$ is a
generating set for $G$ it follows that $G = AH$. The fact that $A
\cap H = 1$ and therefore $G = A \ltimes H$ will follow from the
fact that $H$ is torsion free.

The first level decompositions of $x^{\pm1}$ and $y^{\pm1}$ and the
second level decompositions of $x$ and $y$ are
\begin{align*}
 &x = \s(\s,x) \\
 &y = cabc = \s aaba \s = \s ba \s = x^\s = \s(x,\s) \\
 &x^{-1} = \s(x^{-1},\s) \\
 &y^{-1} = \s(\s,x^{-1})\\
 &x = \s(\s(1,1),\s(\s,x))= \mu(1,1,\s,x) \\
 &y = x^\s = \mu(\s,x,1,1),
\end{align*}
where $\mu=\s(\s,\s)$ permutes the first two levels of the tree as
$00 \leftrightarrow 11, \ 10 \leftrightarrow 01$, which we encode as
the permutation $\mu=(03)(12)$. For a word $w$ over
$\{x^{\pm1},\s\}$ let $\#_x(w)$ and $\#_\s(w)$ denote the total
number of appearances of $x$ and $x^{-1}$ and the number of
appearances of $\s$ in $w$, respectively. Note that a tree
automorphisms represented by a word $w$ over $\{x^{\pm1},\s\}$
cannot be trivial unless both $\#_x(w)$ and $\#_\s(w)$ are even (we
use this claim freely in the further considerations). This is
because $x$ and $x^{-1}$ act as the permutation $(03)(12)$ on the
second level, $\s$ acts as $(02)(13)$, these permutations have order
2, commute, and their product is $(01)(23)$, which is not trivial.

Let $g$ be an element of $H$ that can be written as $g=z_1 z_2 \dots
z_n$, for some $z_i \in \{x^{\pm1},y^{\pm1}\}$, $i=1,\dots,n$.

We claim that if $n$ is odd, the element $g$ cannot have order 2.
Assume otherwise. For $z$ in $\{x^{\pm1},y^{\pm1}\}$ denote $z'=\s
z$. Thus $x'=(\s,x), \ y'=(x,\s)$ and so on. We have
\[ g^2= (z_1 z_2 \dots z_n)^2 = (z_1')^\s z_2' (z_3')^\s z_4' \dots (z_n')^\s z_1'
 (z_2')^\s\dots z_n' =(w_0,w_1),\]
where the words $w_i$ over $\{x^{\pm1},\s\}$ have the property that
\begin{equation}\label{odd}
 \#_x(w_i) = \#_\s(w_i) = n,
\end{equation}
for $i=1,2$. This is because exactly one of $z_i'$ and $(z_i')^\s$
contributes $x^{\pm1}$ to $w_0$ and $\s$ to $w_1$, respectively,
while the other contributes the same letters to $w_1$ and $w_0$,
respectively. But now $g^2$ cannot be 1, since~(\ref{odd}) shows
that neither $w_0$ nor $w_1$ can be 1.

Assume $H$ is not torsion free. Then, since an automorphism of a
binary tree can only have 2-torsion, there exists an element of
order 2 in $H$ and let $g$ be such an element of the shortest
possible length over $\{x,y\}$. Let this length be $n$ and let
$g=z_1 z_2 \dots z_n$, for some $z_i \in \{x^{\pm1},y^{\pm1}\}$,
$i=1,\dots,n$.

Since $n$ must be even, we have
\[ g = z_1 z_2 \dots z_n = (z_1')^\s z_2' \dots (z_{n-1}')^\s z_n' =(w_0,w_1),\]
where $w_0$ and $w_1$ are words over $\{x^{\pm1},\s\}$ whose orders
in $H$ divide 2 and the order of at least one of them is 2. We have
\begin{equation}\label{even}
 \#_x(w_0) \equiv \#_\s(w_0) \equiv  \#_x(w_1) \equiv \#_\s(w_0) \mod 2.
\end{equation}
Indeed, $\#_x(w_i) \equiv \#_\s(w_i) \mod 2$, because $\#_x(w_i) +
\#_\s(w_i) = n$ is even. On the other hand, whenever $z_i'$ or
$(z_i')^\s$ contributes $x^{\pm1}$ or $\s$ to $w_0$, respectively,
it also contributes $\s$ or $x^{\pm1}$ to $w_1$, respectively. Thus
$\#_x(w_0) = \#_\s(w_1)$ and $\#_\s(w_0)=\#_x(w_1)$.

If all the numbers in~(\ref{even}) are even we have that $w_0$ and
$w_1$ represent elements in $H$ (due to the fact that $x^\s=y$) and
can be rewritten as words over $\{x^{\pm1},y^{\pm1}\}$ of lengths at
most $\#_x(w_0) = n - \#_\s(w_0)$ and $\#_x(w_1) = n - \#_\s(w_1)$,
respectively. Either both of these lengths are shorter than $n$ and
therefore none of them can represent an element of order 2 in $H$ or
one of the words $w_i$ is a power of $x$ and the other is trivial.
However, $x$ has infinite order, thus again showing that
$g=(w_0,w_1)$ cannot have order 2.

Finally, assume that all the numbers in~(\ref{even}) are odd. Then,
for $i=1,2$, $w_i$ can be rewritten as $\s u_i$, where $u_i$ are
words over $\{x^{\pm1},y^{\pm1}\}$ of odd length. Let $w_0=\s u_0=\s
t_1 \dots t_m$, where $m$ is odd, and $t_j$ are letters in
$\{x^{\pm1},y^{\pm1}\}$, $j=1,\dots,m$. Then
\[ w_0 = t_1'(t_2')^\s \dots (t_{m-1}')^\s t_m' = (w_{00},w_{01}), \]
where $w_{00}$ and $w_{01}$ are words of odd length $m$ over
$\{x^{\pm1},\s\}$, exactly one of which has even number of $\s$'s.
That word can be rewritten as a word over $\{x^{\pm1},y^{\pm1}\}$ of
odd length. However, the element in $H$ represented by such a word
cannot have order 2, as proved above.
\end{proof}

\begin{lemma}\label{not-metabelian}
The group $H$ is not metabelian.
\end{lemma}

\begin{proof}
Assume that $H$ is metabelian. Consider $[x,y] = x^{-1}y^{-1}xy =
(x^2,x^{-2})$ and $[x,y^{-1}] = x^{-1}yxy^{-1} = (y^2,x^{-2})$. We
have
\begin{align*}
 [[x,y],[x,y^{-1}]] = 1 \qquad &\Rightarrow\\
 [ (x^2,x^{-2}), (y^2,x^{-2})] = 1 \qquad &\Rightarrow \text{(consider the 0 coordinate)}\\
 [x^2,y^2]=1 \qquad &\Rightarrow\\
 [ (x\s, \s x), (\s x, x\s) ] = 1 \qquad &\Rightarrow \text{(consider the 0 coordinate)}\\
 [ x \s, \s x] = 1 \qquad &\Rightarrow\\
 [ (x,\s), (\s,x) ] = 1 \qquad &\Rightarrow \text{(consider the 0 coordinate)}\\
 [ x, \s ] = 1. \qquad &
\end{align*}
The last commuting relation then implies that $x=x^\s=y$, which in
turn implies that $x=\s$, which is impossible since $x$ has infinite
order.
\end{proof}

\begin{lemma}\label{H-branching}
The group $H''$ geometrically contains $H'' \times H'' \times H''
\times H''$, i.e.,
\[ H'' \times H'' \times H'' \times H'' \preceq H''. \]
\end{lemma}

\begin{proof}
The equalities
\begin{align*}
x^2        &= (x, \s, \s, x) \\
y^{-1}x^2y &= (y, x^{-1} \s x, \s, x)
\end{align*}
imply that
\[ H'' \times \langle \s, x^{-1} \s x \rangle'' \times \langle \s
\rangle'' \times \langle x \rangle'' \preceq H''. \]

Since $\langle \s, x^{-1} \s x \rangle$ is dihedral (both generators
have order 2) and $\langle \s \rangle$ and $\langle x \rangle$ are
cyclic, we have
\[ H'' \times 1 \times 1 \times 1 \preceq H''. \]

It follows that $x^{-1} (H'' \times 1 \times 1 \times 1) x = 1
\times 1 \times 1 \times H'' \preceq H''$.

Similarly, the equalities
\begin{align*}
y^2        &= (\s, x, x, \s) \\
x^{-1}y^2x &= (\s, x, y, x^{-1} \s x)
\end{align*}
imply that $ 1 \times 1 \times H'' \times 1 \preceq H''$ and then
$y^{-1}( 1 \times 1 \times H'' \times 1)y = 1 \times H'' \times 1
\times 1 \preceq H''$.
\end{proof}

\begin{corollary}
The group $G$ is a weakly branch group.
\end{corollary}

\begin{proof}
The group $G$, being an infinite self-similar group acting on a
binary three, acts spherically transitively (see Proposition 2
in~\cite{bondarenko_gkmnss:clas32}). The group $H''$ is normal in
$G$, since it is characteristic in the normal subgroup $H$.
Lemma~\ref{not-metabelian} implies that $H''$ is not trivial and
Lemma~\ref{H-branching} then implies that $G_{783}$ is a weakly
branch group.
\end{proof}

\noindent\textbf{802}: $C_2\times C_2\times C_2$. For $G_{802}$ we
have $a=\sigma(a,a)$, $b=(c,c)$, and $c=(a,a)$.

Straightforward calculations give that the group is abelian and has
the indicated structure (there are 8 elements in all and the
nontrivial ones have order 2).

A quick way to see that the group $G_{802}$ is finite is to first
note that the automaton consisting of the states $a$ and $c$
generates the Klein Viergruppe $C_2\times C_2$
(see~Theorem~\ref{thm:class22} and its proof
in~\cite{gns00:automata}) and all sections of $b$ at first level
belong to this automaton. More generally, the following proposition
holds (see, for example,~\cite{bondarenko_s:sushch}).

\begin{prop}
Let $G$ be a self-similar group of tree automorphisms of $X^*$, let
$S$ be a finite generating set for $G$, and let $F$ be a finite set
of automorphisms from $\Aut X^{*}$, such that there exists a level
$k$ in the tree with the property that all sections of all
automorphisms in $F$ at this level belong to $G$ (in particular, $F$
could be a set of finitary automorphisms). Then
\[
\gamma_G(n)\precsim\gamma_{\langle
G,F\rangle}(n)\precsim\bigl(\gamma_G(n)\bigr)^{|X|^k}.
\]
where $\gamma_G(n)$ and $\gamma_{\langle G,F\rangle}(n)$ are the
growth functions of the groups $G$ and $\langle G,F\rangle$ with
respect to the generating sets $S$ and $S\cup F$ respectively.
\end{prop}

\noindent\textbf{843}: For $G_{843}$, we have $a=\sigma(c,b)$,
$b=(a,b)$, and $c=(b,a)$.

The element $c^{-1}a = \s(a^{-1}c,1)$ acts spherically transitively
on the tree (it is conjugate to the adding machine) and therefore
has infinite order. Also, $c^{-1}a$ has infinite order because
$(c^{-1}a)^2=(a^{-1}c,a^{-1}c)$. Since $(c^{-1}ab^{-1}a)^2$ fixes
the vertex $000$ and its section at this vertex equals to $c^{-1}a$,
we get that $c^{-1}ab^{-1}a$ also has infinite order. But the
element $c^{-1}ab^{-1}a$ has itself as a section at vertex $10$.
Hence, $G_{843}$ is not contracting.

\medskip

\noindent\textbf{846}: $C_2\ast C_2\ast C_2$. For $G_{846}$ we have
$a=\sigma(c,c)$, $b=(a,b)$, and $c=(b,a)$.

The automaton 846 is sometimes called Bellaterra automaton, because
it was investigated during the Advanced Course on Automata Groups,
held in Bellaterra, Spain in 2004. A proof was given
in~\cite{nekrash:self-similar}, and here we present the original
proof, which uses action of the automaton dual to the automaton 846.

A \emph{dual} automaton is the automaton that is obtained after
swapping the alphabet and the set of states, as well as the
transition and output functions, in a given automaton. Namely, if
$\mathcal{A} = (Q, X, \tau, \rho)$, then its dual automaton is
$\mathcal{A}^* = (X, Q, \overline{\rho}, \overline{\tau})$, where
$\overline{\rho}(x, q) = \rho(q, x)$, $\overline{\tau}(x, q) =
\tau(q, x)$. The dual automaton represents transitions in the
automata $\mathcal{A}^n$, $n \geq 0$, when they consume input
letters from $X$: the image of the action of $x \in X$ on a word
$q_1q_2\ldots q_n \in Q^n$ is the state $q_1'q_2'\ldots q_n'$ of
$\mathcal{A}^n$ reached after processing the input letter $x$
starting from the initial state $q_1q_2\ldots q_n$. In particular,
the orbits of the semigroup generated by the automaton
$\mathcal{A}^*$ represent the states reachable from a given state in
the automata $\mathcal{A}^n$, $n \geq 0$.

Notion of dual automaton appeared in \cite{aleshin:free} where it
was used to prove that the group generated by two initial automata
is free (unfortunately the proof was not complete, and a complete
proof can be found in~\cite{vorobets:series_free}). A particularly
interesting case is when both $\mathcal{A}^*$ and
$(\mathcal{A}^{-1})^*$ are invertible, such an automaton
$\mathcal{A}$ is called \emph{bi-reversible}. Examples of such
automata are automaton $2240$ generating a free group with three
generators~\cite{vorobets:aleshin}, automaton 846 (under
consideration), and various automata constructed in
\cite{gl_mo:compl}, generating free groups of various ranks.

If the automaton $\mathcal{A}^*$ is invertible then its action on
the set $Q^*$ is in fact a group action, which coincides with the
action of the tree automorphisms group generated by $\mathcal{A}^*$
and we can apply existing techniques to study the action of
$\mathcal{A}^*$ on sets $Q^n$.

Consider the group $G_{846}$. Its generators $a$, $b$, and $c$ have
order $2$, therefore to prove that $G_{846}$ is isomorphic to $C_2
\ast C_2 \ast C_2$ we need to show that no word in $w \in W_n$, $n
\geq 1$, is trivial in $G_{846}$, where $W_n \subset \{a, b,
c\}^{*}$ is the set of all words of length $n$ which do not contain
squares of the letters $a$, $b$, and $c$. For any $n > 0$, the set
of words in $W_n$ that are not trivial in $G_{846}$ is nonempty. For
instance, such a word is $w = abcbcb\cdots$ (this is because $b$ and
$c$ act trivially on the level 1 of the tree, while $a$ permutes the
two vertices at level 1). If the state $w$ in the $n$-th power of
the automaton 846 is reachable from some state $w' \in W_n$ then the
element $w'$ is not trivial in $G_{846}$. Therefore it is sufficient
to show that the dual automaton acts transitively on the sets $W_n$,
$n \geq 1$.

The dual automaton to $846$ is the invertible automaton defined by
\begin{equation}
\begin{array}{rcr}\label{dual-bellaterra}
A & = & (acb)(B, A, A), \\
B & = & (ac)(A, B, B).
\end{array}\end{equation}

The group $D =\langle A, B \rangle$ does not act spherically
transitively on the ternary tree over the alphabet $\{a,b,c\}$ (for
instance the orbit of $a^2$ under action of $D$ consists only of
$a^2$, $b^2$, and $c^2$). Let us prove that $D$ acts transitively on
the sets $W_n$. The set $T' = \bigcup_{n \ge 0}W_n$ of all words
that do not contain squares of the letters in $\{a,b,c\}$ is a
subtree of the ternary tree $\{a, b, c\}^*$ and this subtree is
invariant under the action of $D$. The root of $T'$ is connected to
the vertices $a$, $b$ and $c$, which, in turn, are roots of $3$
binary trees. The action of $D$ on this subtree can be derived from
the relations~\eqref{dual-bellaterra}. The generators $A$ and $B$
act as follows:
\begin{equation}
\begin{array}{rcr}
A & = & (acb)(B_a, A_b, A_c) \\
B & = & (ac)(A_a, B_b, B_c)
\end{array}
\end{equation}
where $A_a, A_b, A_c, B_a, B_b, B_c$ are the automorphisms of the
binary trees hanging down from the vertices $a$, $b$ and $c$. All
these binary trees can be naturally identified with the tree
$\{0,1\}^*$, where the action of $A_a, A_b,\ldots,B_c$ is defined by
\begin{equation}\begin{array}{rcr}
A_a & = & (A_b, A_c), \\
A_b & = & \sigma(B_a, A_c), \\
A_c & = & \sigma(B_a, A_b), \\
B_a & = & \sigma(B_b, B_c), \\
B_b & = & \sigma(A_a, B_c), \\
B_c & = & \sigma(A_a, B_b). \\
\end{array}
\end{equation}

The algorithm deciding the spherical transitivity of an element
acting on the binary tree (see Algorithm 1
in~\cite{bondarenko_gkmnss:clas32}) shows that the element $B_b$
acts level transitively. Thus $D$ acts transitively on the levels of
$T'$. Indeed, $D$ acts transitively on the first level, $B$
stabilizes vertex $b$ and its section at $b$ is $B_b$.

\medskip

\noindent\textbf{849}: For $G_{849}$ we have $a=\sigma(c,a)$,
$c=(1,a)$ and $b=1$.

The element $a^2c=(ac,ca^2)$ is nontrivial because its section at
vertex $0$ is $ac$, which acts nontrivially on the first level. The
section of $(a^2c)^2$ at vertex $00$ coincides with $a^2c$, hence
this element has infinite order. On the other hand, the section of
$a^2c$ at vertex $100$ coincides with $a^2c$, which shows that
$G_{849}$ is not contracting.

The group $G_{849}$ is regular weakly branch over $G_{849}'$ because
it is self-replicating and $[a^{-1},c]\cdot[c,a]=([a,c],1)$.

Conjugating the generators of $G_{849}$ by the automorphism
$\mu=\s(\mu,c^{-1}\mu)$ yields the recursion
\[ x = \s(yx,1), \qquad y = (x,1), \]
where $x = a^\mu$ and $y = c^\mu$. Since
\[ x = \s(yx,1), \qquad yx = \s(yx,x), \]
and the last recursion defines the automaton 2852, we see that
$G_{849} \cong G_{2852}$. See $G_{2852}$ for further information.

\medskip

\noindent\textbf{875}.  For $G_{875}$ we have $a=\sigma(b,a)$,
$b=(b,c)$, and $c=(b,a)$.

Let us prove that $a$ has infinite order. This will be accomplished
by proving that the orbit of $110^\infty\in X^\omega$ under the
action of $a$ is infinite.

For every $w\in X^*$,
$a(w0^\infty)=a(w0)a|_{w0}(0^\infty)=a(w0)b(0^\infty)=a(w0)0^\infty$.
Therefore all points in the orbit of $110^\infty$ under the action
of $a$ end in $0^\infty$. For a word $w$ ending in $0^\infty$ define
the \emph{head} of $w$ to be shortest word $u$ such that $w =
u0^\infty$. The length of the head of any infinite word ending in
$0^\infty$ cannot decrease under the action of $a$. Indeed, for
every $w\in X^*$,
\begin{align*}
 a(w010^\infty) &=a(w0)a|_{w0}(10^\infty)=a(w0)b(10^\infty)=a(w0)10^\infty
 \quad\text{and} \\
 a(w110^\infty) &=a(w11)a|_{w11}(0^\infty)=a(w11)a(0^\infty)=a(w11)10^\infty.
\end{align*}
On the other hand, the length of the head along the orbit cannot
stabilize, because in this case the orbit would be finite and we
must have $a^k(110^\infty)=110^\infty$, for some $k\geq1$. But this
is impossible since $a(110^\infty)=0010^\infty$ and the length of
the head does not decrease. Thus the orbit is infinite and $a$ has
infinite order.

Since $c=(b,a)$ and $b=(b,c)$, the elements $b$ and $c$ also have
infinite order. Finally, since $b=(b,c)$ and $b$ has infinite order,
$G_{875}$ is not contracting.

\medskip

\noindent\textbf{891}: $C_2\ltimes Lamplighter$. For $G_{891}$ we
have $a=\sigma(c,c)$, $b=(c,c)$, and $c=(b,a)$.

All generators have order $2$.

Consider the subgroup $H=\langle x, y \rangle$, where $x = ac$ and
$y = cb$. We will prove that $H$ is isomorphic to the Lamplighter
group $\Z \ltimes C_2$. We have
\[
 x = ac = \s(cb,ca) = \s(y,x^{-1}), \qquad
 y = cb = (bc,ac) = (y^{-1},x).
\]
Note that $xy=\s$.

Consider the elements $s_n=\sigma^{y^n}=y^{-n}xy^{n+1}$, $n\in\Z$.
We have, for $n >0$, $s_0s_1\cdots s_{n-1}=x^ny^n$ and
$s_{-n}s_{-n+1}\cdots s_{-1}=y^nx^n$. Further, $s_n=y^{-n}\sigma
y^n=\sigma(x^{-n}y^{-n},y^{n}x^{n})$ and $s_{-n}=y^{n}\sigma
y^{-n}=\sigma(x^ny^n,y^{-n}x^{-n})$. Hence,
\[s_n=\sigma(s_{-1}s_{-2}\cdots s_{-n}, s_{-n}\cdots s_{-2}s_{-1})\]
and
\[s_{-n}=\sigma(s_{0}s_{1}\cdots s_{n-1}, s_{n-1}\cdots s_{1}s_{0}).\]

Now by induction we get that the depth of $s_n$ is $2n+1$ for $n\geq
0$ and the depth of $s_{-n}$ is $2n$ for $n>0$ (by \emph{depth} of a
finitary element we mean the lowest level at which all sections of
the element are trivial). Therefore, all $s_i$ are different,
commute (because for each $i$ and each level $m$ all sections of
$s_i$ at level $m$ are equal), and have order $2$ (they are all
conjugate to $\s$). Note that this proves that $y$ has infinite
order and that $H = \langle x,y \rangle = \langle y, \sigma \rangle
\cong \Z \ltimes C_2$.

The subgroup $H$ has index at most 2 in $G_{891}$ because the
generators of $G_{891}$ are of order $2$. In fact, it has index
exactly 2, since there is no relation of odd length in $G_{891}$.
Indeed, let $w$ be a word over $\{a,b,c\}$ representing the identity
in $G_{891}$. Then $\#_a(w)$ must be even (since $a$ is the only
generator that acts nontrivially on the first level). Further,
$\#_c(w)$ must also be even. This is because $w$ can be decomposed
as $w=(w_0,w_1)$ where $w_0$ and $w_1$ are words over $\{a,b,c\}$
such that $\#_a(w_0) + \#_a(w_1) = \#_c(w)$ (only the generator $c$
contributes $a$ to the next level in the decomposition). Further,
$\#_c(w_0) = \#_c(w_1) = \#_a(w) + \#_b(w)$, since both $a$ and $b$
contribute $c$ to both coordinates on the next level. Thus if
$\#_b(w)$ were odd (and we already know that $\#_a(w)$ is even) then
none of the words $w_0$ and $w_1$ would represent the identity. Note
that we proved by this that every word $w$ representing the identity
in $G_{891}$ must have even number of occurrences of $a$ $b$ and
$c$, which shows not only that $c$ does not belong to $H$, but also
that the abelianization of $G_{891}$ is $C_2 \times C_2 \times C_2$
(and the commutator consists precisely of the elements represented
by words $w$ for which $\#_a(w)$, $\#_b(w)$, and $\#_c(w)$ are
even).

Therefore, the group $G_{891}=\langle c,H \rangle$ has the structure
of a semidirect product $G_{891} = C_2 \ltimes ( \Z \wr C_2)$ and
the action of $c$ on $H$ is given by $(x)^c=(ac)^c = ca = x^{-1}$
and $y^c=(cb)^c = bc = y^{-1}$. It follows that $G_{891}$ is
solvable group of exponential growth.

Since $y$ has infinite order, stabilizes the vertex $00$ and has
itself as a section at this vertex, we conclude that $G_{891}$ is
not contracting.

\medskip

\noindent\textbf{2193}. For $G_{2193}$ we have $a=\sigma(c,b)$,
$b=\sigma(a,a)$, and $c=(a,a)$.

For $x=ca^{-1}$ and $y=ab^{-1}$, we have
$x=\s(ab^{-1},ac^{-1})=\s(y,x^{-1})$ and
$y=(ba^{-1},ca^{-1})=(y^{-1},x)$. It is proved above (see
$G_{891}$), that $\langle x,y\rangle$ is not contracting and is
isomorphic to the Lamplighter group. Thus $G_{2193}$ is neither
torsion, nor contracting and has exponential growth.

\medskip

\noindent\textbf{2280}.  For $G_{2280}$ we have $a=\sigma(c,a)$,
$b=\sigma(b,a)$, and $c=(b,a)$.

Let us prove that $a$ has infinite order. We will prove that the
orbit of $10^\infty$ under iterations of $a^2$ is infinite. Recall
the definition of a head for a word over $\{0,1\}$ ending in
$0^\infty$ (see $G_{875}$). For every word $w \in X^*$, we have
$a^2(w10^\infty)=a^2(w)\cdot a^2|_w(10^\infty)$. The automorphism
$a^2$ has 6 distinct sections defined by
\begin{alignat*}{5}
 a^2 &= (ac,ca), \qquad &&ac &&= \s (cb,a^2), \qquad &&ca &&=\s(ac,ba)
 \\
 cb &= \s(ab,ba), &&ba &&= (ac,ba),  &&ab &&=(ab,ca).
\end{alignat*}
Since
\begin{align*}
 a^2(10^\infty) &= ab(10^\infty) = 1110^\infty, \\
 ac(10^\infty) &= ca(10^\infty) = cb(10^\infty) = 0010^\infty, \text{ and}\\
 ba(10^\infty) &= 10110^\infty
\end{align*}
we see that the length of the head increases by 3 at each
application of $a^2$ and $a$ has infinite order.

Since $a^2=(ac,ca)$, the element $ca$ also has infinite order. The
element $ab$ has infinite order, because $ab=(ab,ca)$. Therefore,
$G_{2280}$ is not contracting.

\medskip

\noindent\textbf{2294}: Baumslag-Solitar group $BS(1,-3)$. For
$G_{2294}$ we have $a=\sigma(b,c)$, $b=\sigma(c,a)$, and $c=(b,a)$.

The automaton satisfies the conditions of Proposition 1
in~\cite{bondarenko_gkmnss:clas32}, which implies that $cb$ has
infinite order. Since $a^2=(cb,bc)$, the generator $a$ also has
infinite order.

Let $\mu = ca^{-1}$. We have
$\mu=ca^{-1}=\sigma(ac^{-1},1)=\sigma(\mu^{-1},1)$, which shows that
$\mu$ has infinite order (it acts spherically transitively on the
binary tree). Since $bc^{-1} = \s(cb^{-1},1) =
\s((bc^{-1})^{-1},1)$, we see that $bc^{-1} = \mu = ca^{-1}$ and
therefore $G_{2294} = \langle \mu, a \rangle$.

We claim that $a\mu a^{-1}=\mu^{-3}$ in $G_{2294}$. We have
\[
 a \mu a^{-1} = aca^{-2} = \s(bc^{-1},cac^{-1}b^{-1}),
  \qquad
 \mu^{-3} = \s(\mu, \mu^2).
\]
Therefore it is sufficient to show that $cac^{-1}b^{-1} =
ca^{-1}ca^{-1}$, i.e., $ac^{-1}b^{-1}ac^{-1}a=1$. This is clear
since $ac^{-1}b^{-1}ac^{-1}a = (ca^{-2}ca^{-1}b,1)$ and
$ca^{-2}ca^{-1}b$ is conjugate of the inverse of
$ac^{-1}b^{-1}ac^{-1}a$. This completes the proof.

For realizations of $BS(1,m)$ for any value of $m$, $m \neq \pm 1$,
see~\cite{bartholdi-s:bs}.

\medskip

\noindent\textbf{2396}: \emph{Boltenkov group}. For $G_{2396}$ we
have $a=\sigma(b,a)$, $b=\sigma(c,b)$, and $c=(c,a)$.

This group was studied by Boltenkov. He proved that it is torsion
free and that the monoid generated by $a$, $b$, and $c$ is free.
Here we provide his proofs and show that the group is not
contracting.

\begin{prop}\label{free}
The monoid generated by $a$, $b$, and $c$ is free.
\end{prop}
\begin{proof}
Suppose there are some relations. Let $w=u$ be a relation for which
$\max(|w|,|u|)$ minimal. Assume that none of the words $w$ and $u$
is empty. Clearly, $w$ and $u$ end in different letters (since
cancelation holds). Then $w=\sigma_w(w_0,w_1)=\sigma_u(u_0,u_1)=u$,
where $\sigma_w,\sigma_u$ are permutations in $\{1,\s\}$. It follows
that $w_0=u_0$ and $w_1=u_1$ are also relations. Consider different
cases.

Let $w$ end in $b$ and $u$ end in $c$. Then $w_0$ and $u_0$ both end
in $c$. This means, by minimality, that $w_0=u_0$ as words. Thus, in
particular, $|u|=|w|$. Since $b \neq c$ in $G_{2396}$ (their actions
differ at level 1) the length of $w$ and $u$ is at least 2. We can
unambiguously recover the second to last letter in $w$ and $u$.
Indeed, the second to last letter in $u_0$ can be only $b$ or $c$
(these are the only possible sections at 0 of the three generators),
while the second to last letter of $w_0$ can be only $a$ or $b$
(these are the only possible sections at 1 of the three generators).
Thus $w_0=u_0=\dots bc$, $w=\dots bb$ , and $u=\dots ac$. As $bb\neq
ac$ in $G_{2396}$ (their actions differ at level 1), the length of
$w$ and $u$ is at least 3. Iterating the procedure, we obtain that
$w_0=u_0=b\dots bbc$, and therefore $w=\dots ababb$, $u=\dots
babac$. As $|u|=|w|$, the actions of $w$ and $u$ differ at level 1
and we get a contradiction.

Let $w$ end in $a$ and $u$ end in $b$ or $c$. Then $u_0,w_0$ end in
$b$ and $c$, respectively, and the situation is reduced to the case
above.

Therefore, there are no nontrivial relations in $G_{2396}$ of the
form $w=u$, for nonempty words $w$ and $u$ over $\{a,b,c\}$.

If $u$ is an empty word, then $w_0=1=w_1$, which shows that
$w_0=w_1$ is also a minimal relation, which is impossible because
both $w_0$ and $w_1$ are nonempty.
\end{proof}

For a group word $w$ over $\{a,b,c\}$, define the exponent
$\exp_a(w)$ of $a$ in $w$ to be the sum of the exponents of all the
occurrences of $a$ and $a^{-1}$ in $w$ (define $\exp_b(w)$ and
$\exp_c(w)$ analogously).  Let $\exp(w)=\exp_a(w)+\exp_b(w) +
\exp_c(w)$. Note that if $w=\sigma_w(w_0,w_1)$, then
$\exp(w)=\exp(w_0)=\exp(w_1)$.

\begin{lemma}
If $w=1$ in $G_{2396}$ then $\exp(w)=0$
\end{lemma}
\begin{proof}
Suppose not. Choose a freely reduced group word $w$ over $\{a,b,c\}$
such that $w=1$ in $G_{2396}$, $\exp(w)\neq 0$, and $w$ has minimal
length among all such words. Let $w=(w_0,w_1)$. Then $w_0,w_1$ also
represent $1$ in $G_{2396}$ and $\exp(w_0)=\exp(w_1)=\exp(w)\neq 0$.
Therefore $w_0$ and $w_1$ are not empty and, by minimality, their
length must be equal to the length of $w$. This implies that $w$
does not contain $ac^{-1}$, $bc^{-1}$, $ca^{-1}$, or $cb^{-1}$ as a
subword (this is because $ac^{-1}=\sigma(bc^{-1},1)$ and
$bc^{-1}=\sigma(1,ba^{-1})$).

By the same argument $w_0$ and $w_1$ must not have the above 4 words
as subwords  as well, implying that $w$ does not have
$ab^{-1}=(ab^{-1},bc^{-1})$ or its inverse $ba^{-1}$ as a subword.
Thus, since $w$ is reduced, we have that
$1=w=W_1(a^{-1},b^{-1},c^{-1})W_2(a,b,c)$, and we obtain a relation
between positive words over $\{a,b,c\}$, contradicting
Proposition~\ref{free}.
\end{proof}

Note that if $w=1$ in $G_{2396}$ then $\exp_a(w)$, $\exp_b(w)$ and
$\exp_c(w)$ are even. Indeed, clearly $\exp_a(w)+\exp_b(w)$ must be
even (by looking at action on level 1), which implies by the
previous lemma that $\exp_c(w)$ is even. Further, if $w=(w_0,w_1)$,
then also $\exp_a(w_0)+\exp_b(w_0)$ and $\exp_a(w_1)+\exp_b(w_1)$
are even. We have that
$\exp_a(w)+\exp_b(w)=\exp_b(w_0)+\exp_b(w_1)$,
$\exp_a(w)+\exp_c(w)=\exp_a(w_0)+\exp_a(w_1)$. Therefore $2\exp_a(w)
+ \exp_b(w) + \exp_c(w)$ is even, implying that $\exp_b(w)$ is even.
Since both $\exp_b(w)$ and $\exp_c(w)$ are even, so must be
$\exp_a(w)$.

\begin{prop}
\label{boltenkov_torsion_free} The group $G_{2396}$ is torsion free.
\end{prop}
\begin{proof}
Suppose not. Then $G_{2396}$ has an element $w$ of order $2$. By
considering the sections we can also suppose that $w$ does not
belong to the stabilizer of the first level, i.e.
$w=\sigma(w_0,w_1)$. It follows that $w^2=(w_1w_0,w_0w_1)=1$. We
have that, modulo 2, $\exp_b(w_0w_1) =  \exp_b(w_0) + \exp_b(w_1)
=\exp_a(w)+\exp_b(w)$. Since $\exp_b(w_0w_1)$ is even, we conclude
that $\exp_a(w)+\exp_b(w)$ is even, contradicting the assumption
that $w$ does not stabilize the first level.
\end{proof}

The element $b^{-1}a=(c^{-1}b,b^{-1}a)$ is nontrivial (the $a$
exponent is odd) and therefore has infinite order by
Proposition~\ref{boltenkov_torsion_free}. On the other hand,
$b^{-1}a$ fixes vertex 1 and its section at this vertex is
$b^{-1}a$, implying that $G_{2396}$ is not contracting.

\medskip

\noindent\textbf{2398}: This group was studied by Dahmani
in~\cite{dahmani:non-contr}. It is self-replicating, not
contracting, weakly regular branch group over its commutator
subgroup.

\medskip

\noindent\textbf{2852=$G_{849}$}. For $G_{2852}$ we have
$a=\sigma(b,1)$, $b=\sigma(b,a)$, $c=1$.

Since $a^2 = (b,b)$ and $ab=(b,ba)$, the group $G_{2852}$ is
self-replicating and spherically transitive.

The group $G_{2852}$ is a regular weakly branch group over
$G_{2852}'$. Indeed, $[a^{-1},b]\cdot[b,a]=([a,b],1)$. Spherical
transitivity and the self-replicating property then imply that
$G_{2852}'\times G_{2852}'\preceq G_{2852}'$. Since $G_{2852}$ is
not abelian ($[b,a] = (b^{-1}ab,a^{-1}) \neq 1$, because $a \neq
1$), $G_{2852}'$ is nontrivial and $G_{2852}$ is a regular weakly
branch over $G_{2852}'$.

We have $b^2=(ab,ba)$, $ba=(ab,b)$, and $ab=(b,ba)$. Thus $b^2$
fixes the vertex $00$ and has $b$ as a section at that vertex. Since
$b$ is nontrivial, this implies that $b$ has infinite order, and
therefore so does $ab$. On the other hand, $ab$ fixes vertex $10$
and has itself as a section at that vertex, implying that $G_{2852}$
is not contracting.

We claim that the monoid generated by $a$ and $b$ is free (and
therefore the group has exponential growth).

Indeed, let $w$ be a nonempty word in $\{a,b\}^*$. If $w=1$ in
$G_{2852}$, then $w$ contains both $a$ and $b$, because they both
have infinite order ($a$ has infinite order since $a^2=(b,b)$).
Suppose the length of $w$ is minimal among all nonempty words over
$\{a,b\}$ representing the identity element in $G_{2851}$. Then one
of the sections of $w$ will be shorter than $w$ (since $a|_1$ is
trivial), nonempty (since $b|_0$ and $b|_1$ are nontrivial), and
will represent the identity in $G_{2851}$, which contradicts the
minimality assumption. Thus $w\neq 1$ in $G_{2851}$, for any
nonempty word in $\{a,b\}^*$.

Now take any nonempty words $w,v\in\{a,b\}^*$ with minimal sum
$|w|+|v|$ such that $w=v$ in $G_{2852}$. Let $w=\s_w(w_0,w_1)$ and
$u=\s_u(u_0,u_1)$, where $\s_w,\s_w \in \{1,\s\}$. Suppose $w$ ends
in $a$ and $v$ ends in $b$. Then $w_1=v_1$ in $G$. If $w=a$ then
$w_1=1$ in $G$ and the word $v_1$ is nontrivial, which is
impossible. If $w$ contains more than one letter then $w_1$ ends in
$b$, $v_1$ ends in $a$, and $|w_1|<|w|$ and $|v_1| \leq |v|$. Thus
we have a shorter relation $w_1=v_1$, contradicting the minimality
assumption.

See $G_{849}$ for the isomorphism between $G_{2852}$ and $G_{849}$.

Conjugating the generators of $G_{2852}$ by the automorphism
$\mu=\s(b\mu,\mu)$ yields the recursion
\[ x = \s(y,1), \qquad y = \s(xy,1), \]
where $x = a^\mu$ and $y = b^\mu$. Since
\[ y = \s(xy,1), \qquad xy = (xy,y), \]
and the last recursion defines the automaton 933, we see that
$G_{2852} \cong G_{933}$.

\bibliographystyle{amsalpha}

\newcommand{\etalchar}[1]{$^{#1}$}
\def\cprime{$'$} \def\cprime{$'$} \def\cprime{$'$} \def\cprime{$'$}
  \def\cprime{$'$}
\providecommand{\bysame}{\leavevmode\hbox
to3em{\hrulefill}\thinspace}
\providecommand{\MR}{\relax\ifhmode\unskip\space\fi MR }
\providecommand{\MRhref}[2]{%
  \href{http://www.ams.org/mathscinet-getitem?mr=#1}{#2}
} \providecommand{\href}[2]{#2}

\end{document}